\pdfoutput=1
\documentclass[11pt,a4paper]{article}
\usepackage[T1]{fontenc}
\usepackage[utf8]{inputenc}
\usepackage{lmodern}
\usepackage[margin=25mm]{geometry}
\usepackage{amsmath,amssymb,amsthm,mathtools,bm}
\usepackage{graphicx,booktabs,microtype,xcolor}
\usepackage[colorlinks=true,linkcolor=blue!45!black,citecolor=blue!45!black,urlcolor=blue!45!black]{hyperref}
\usepackage{enumitem}
\usepackage{placeins,longtable,array,fancyhdr}
\definecolor{deepteal}{HTML}{006D77}
\hypersetup{linkcolor=deepteal,citecolor=deepteal,urlcolor=deepteal,
pdftitle={Boundary Actuation of a Naturally Stable Chain Fountain},pdfauthor={Latchezar Tomov},
pdfkeywords={chain fountain, Mould effect, boundary actuation, nonlinear tracking, state constraints, robust control, Lyapunov stability, stochastic viability, backstepping, zero-gravity circulation, Lean 4 formalization}}
\newcommand{\lean}[1]{{\footnotesize\path{#1}}}

\setlist{itemsep=2pt,topsep=5pt}
\numberwithin{equation}{section}
\newtheorem{theorem}{Theorem}[section]
\newtheorem{proposition}[theorem]{Proposition}
\newtheorem{corollary}[theorem]{Corollary}
\theoremstyle{definition}
\newtheorem{assumption}[theorem]{Assumption}

\newcommand{\R}{\mathbb R}
\newcommand{\D}{\mathcal D_+}
\newcommand{\dd}{\mathrm d}
\newcommand{\norm}[1]{\left\lVert#1\right\rVert}
\newcommand{\amin}{\alpha_-}
\newcommand{\amax}{\alpha_+}
\newcommand{\un}{u_{\mathrm n}}
\newcommand{\ur}{u_{\mathrm r}}
\newcommand{\Vbar}{\overline V}
\newcommand{\JitterReduction}{87.99}
\newcommand{\RefinementError}{1.15e-14}
\newcommand{\MinSimTension}{19.5392}
\newcommand{\MinSimCrown}{4.7497}
\newcommand{\MaxSimV}{0.000163}
\newcommand{\MaxSimInput}{0.8114}

\title{Boundary Actuation of a Naturally Stable Chain Fountain:\\Mechanics, Physically Admissible Tracking,\\Stochastic Limits and Zero-Gravity Circulation}
\author{Latchezar Tomov\\\small Department of Informatics, New Bulgarian University, Sofia, Bulgaria\\
\small \href{mailto:lptomov@nbu.bg}{lptomov@nbu.bg}\quad
\href{https://orcid.org/0000-0003-1902-6473}{ORCID: 0000-0003-1902-6473}}
\date{Expanded preprint with mathematical appendix and zero-gravity extension\\14 November 2026}
\begin{document}
\maketitle
\begin{abstract}
This paper was inspired by a post on Twitter and by the eight months the author spent at the AASS Centre at \"Orebro University in 2007. This phenomenon is interesting not only to a physicist, but also to someone like me working in control theory, because of its potential relevance to stabilizing drones and similar flying machines subject to random velocity fluctuations while hovering. A chain fountain can settle without an engineered controller, yet its natural height is selected by transport and boundary forces. We ask how an ideal boundary traction force can sustain a prescribed height or motion. Within an established three-state model, we connect feedback design to force and energy balances and derive sufficient conditions for nominal tracking that preserve the model's mechanical constraints. We also bound the deviations caused by bounded pickup uncertainty. The resulting framework explains how feedback can shape a naturally regulated process and identifies physical limits on the motions it can sustain: stabilizing a mathematical trajectory is useful only if the corresponding motion remains mechanically admissible. This is a preliminary step towards stabilization of the chain in three dimensions.
\end{abstract}

\medskip
\noindent\textbf{Public resources:}
\href{https://luchesart.github.io/Chain-fountain/}{Interactive simulator};
\href{https://github.com/luchesart/Chain-fountain}{project repository};
\href{https://github.com/luchesart/Chain-fountain/blob/main/Chain_Fountain_Lean_GitHub.zip}{Lean source archive};
\href{https://luchesart.github.io/Chain-fountain/Chain_Fountain_Zero_Gravity_Proof_Map.html}{interactive proof map}.

\clearpage
\tableofcontents
\clearpage
\section{Introduction}
When a chain flows from an elevated container to the floor, part of it may rise above the container before turning downward. The question is not merely why a moving chain can form an arch. It is also how the arch reaches a persistent height, which quantities regulate that height, and how an imposed force changes the balance. These questions lead to a control formulation in which the controller must respect the mechanical assumptions of the model.

We start from Biggins's transient equations for the rising leg, falling leg and moving apex~\cite{biggins}. An ideal signed force at pickup becomes the control input. The objective is to determine when that force can move the fountain toward a prescribed height while the chain continues to satisfy the transport assumptions from which the equations were derived.

Two distinctions organize the paper. First, convergence of a height error does not imply that the underlying chain continues to move in its assumed direction. A cancellation controller can track a mathematically prescribed curve while leaving the physical domain. Second, a bounded random coefficient is different from white noise acting directly on velocity. A bound that is valid for every bounded coefficient realization cannot be transferred to an It\^o diffusion without examining the diffusion at the safety boundary.

\subsection*{Mechanical and control-theoretic background}
The mechanical and control questions meet at the boundary of the moving chain. The force balance along a freely moving segment describes how tension, curvature and gravity are related, but the complete motion also depends on how chain enters and leaves that segment. Biggins and Warner's steady analysis makes these boundary reactions part of the height selection mechanism, and Biggins's transient model follows the evolution toward that operating point~\cite{bigginswarner,biggins}. We use the latter model because it separates apex motion from the rising and falling transport speeds. That separation allows an imposed force to change the height while leaving an internal speed whose evolution must still be examined.

The choice of boundary law requires particular care. A continuum description can represent rapid pickup and deposition through jump balances and dissipation, as in Virga's treatment of dissipative shocks~\cite{virga2014}. At the scale of individual links, the contact geometry determines how the chain is brought into motion. Pantaleone compared a pickup model with measured speeds and heights for several chains~\cite{pantaleone2017}, while Martins's experiments and the simulations and experiments of Flekk\o y, Moura and M\aa l\o y show why chain arrangement, flexibility and contact with the underlying packing matter~\cite{martins2019,flekkoy2018,flekkoy2019}. Pfeiffer and Mayet also studied the stationary fountain through both continuum and multibody descriptions~\cite{pfeiffer2017}. These works explain why our pickup coefficient is an effective parameter. Allowing a bounded variation of that coefficient tests robustness within the chosen model, but does not replace a model of the changing contacts that produce the variation.

There is also a direct mechanical precedent for deliberately selecting a fountain. Mayet and Pfeiffer proposed a system of guiding pulleys and weighted carts that produces prescribed stationary boundary forces~\cite{mayet2020}. Their construction starts with a desired fountain geometry, calculates its required forces and chooses the apparatus parameters accordingly. It therefore already connects a selected fountain to manipulable mechanical conditions. Our question concerns the subsequent dynamics: given an ideal boundary force servo and a changing height reference, which initial states and force limits permit the reference to be followed? The difference is between constructing a stationary force balance and proving the behavior of a feedback system during an excursion. Both depend on the physical meaning of the boundary forces.

The relevant control literature extends well beyond fountains. Axially moving strings and belts have long provided models for transport with coupled motion and vibration. Nguyen and Hong regulate the transport velocity while suppressing longitudinal and transverse vibrations through a boundary force and torques at two drive rollers~\cite{nguyen2012}. Constraints have also been incorporated into this setting. Liu, Zhao and Guo combine adaptive backstepping with an auxiliary system for actuator saturation~\cite{liu2016}, and Zhao and coauthors use a barrier Lyapunov function to enforce a boundary tension constraint in the presence of an unknown disturbance~\cite{zhao2017}. These results establish important precedents for velocity tracking, robust boundary action and constrained string control. Their distributed models and actuator equations address degrees of freedom that the three-state fountain approximation does not contain. Conversely, the fountain model has its own pickup law, moving crown and positive-transport domain, which require separate inequalities.

The distinction between a reduced model and a spatial model matters in the stability argument itself. Zhao and Rahn showed that, for a pinned axially moving string, an inappropriate use of the material derivative in a Lyapunov calculation and omission of a coupled domain can lead to incorrect conclusions~\cite{zhaorahn2006}. Here the tracking estimate is obtained by differentiating a function of the actual error coordinates along the finite-dimensional closed-loop equations. The mechanical energy balance separately accounts for transport across the open boundaries. When spatial deformation is introduced, its energy fluxes, boundary conditions and coupling to the actuator must be included in a new analysis. Stability of the reduced error equations supplies no estimate for those additional modes.

Within this setting, backstepping and control barrier functions provide established design tools~\cite{krstic,ames}. The purpose of the present analysis is to connect their conclusions to the particular mechanical domain. For example, decreasing height error can coexist with a falling speed approaching zero. A bound on the error must therefore be combined with a comparison argument for the internal speed. The resulting height and speed bounds also keep the equations away from their singular denominators, which makes continuation of the locally constructed solution possible. Only after this step do the tension and force estimates yield an all-time admissible tracking statement under the stated hypotheses. This sequence explains the role of the sufficient conditions and gives a concrete basis for comparison with other control designs for the same model.

The zero-gravity extension changes the mechanical problem again. Removing gravity from open transport preserves pickup and deposition, whereas a closed circulating string has periodic material geometry and no such exchange with stationary piles. Motions along a fixed curve with constant circulation and tension belong to classical moving-string mechanics~\cite{hanna2013}. Section~\ref{sec:zerog} places a speed regulator on that exact family by specifying a spatially distributed tangential force. Its scalar tracking estimate describes changes of circulation within the family. Shape perturbations and actuation through a localized drive lead to different equations, so comparison with the spatial control literature must be made at the level of the plant and the admissible forces as well as the selected output.

\subsection{Why control a process that already settles?}\label{sec:whycontrol}
The unforced fountain is an instructive example of a stable process whose operating point is not freely assignable. Its height is selected jointly by gravity, transport, pickup and deposition. In the adopted model the natural height is
\[
h^*=\frac{\alpha_0H}{1-\alpha_0-\beta}.
\]
Local stability means that sufficiently small perturbations return to this particular operating point. It does not say that a different height can be maintained, that the same transient response is obtained when the packing changes, or that an arbitrarily rapid height request preserves forward transport. Those are distinct control questions.

\textbf{There is no designed regulator in the unforced experiment.} The restoring action arises from the coupled mechanical balances. When the rising and falling speeds cease to match, the apex moves, changing the leg lengths and hence their accelerations. At the same time, speed-dependent boundary forces and losses affect the transport. The sign of any one of these effects does not establish stability, so Section~\ref{sec:spectrum} examines their complete local coupling through the characteristic polynomial and a nonlinear Lyapunov argument.

\textbf{The regulator introduced here is a feedback law commanding a boundary force.} Its manipulated variable is $u$, the force per unit linear density exerted downward on the moving chain near pickup. From the state $(h,v_1,w)$ and the desired height together with its first two derivatives, the controller computes a force that makes the height-error dynamics dissipative. A positive command brakes the chain, while a negative command assists its motion. The effective pickup coefficient $\alpha$ is a model parameter, not a control input. An experiment based on this formulation would require a force-controlled traction device whose influence on the pickup law is characterized separately.

This formulation allows us to ask which steady heights are compatible with a given force range and which moving references can be followed while both legs continue transporting chain in the assumed directions. It also allows us to quantify how much disturbance rejection is possible before the required force or pickup tension becomes inadmissible. Stability of the uncontrolled plant provides a useful starting point for these questions, but does not answer them.

\begin{figure}[htbp]
\centering\includegraphics[width=\linewidth]{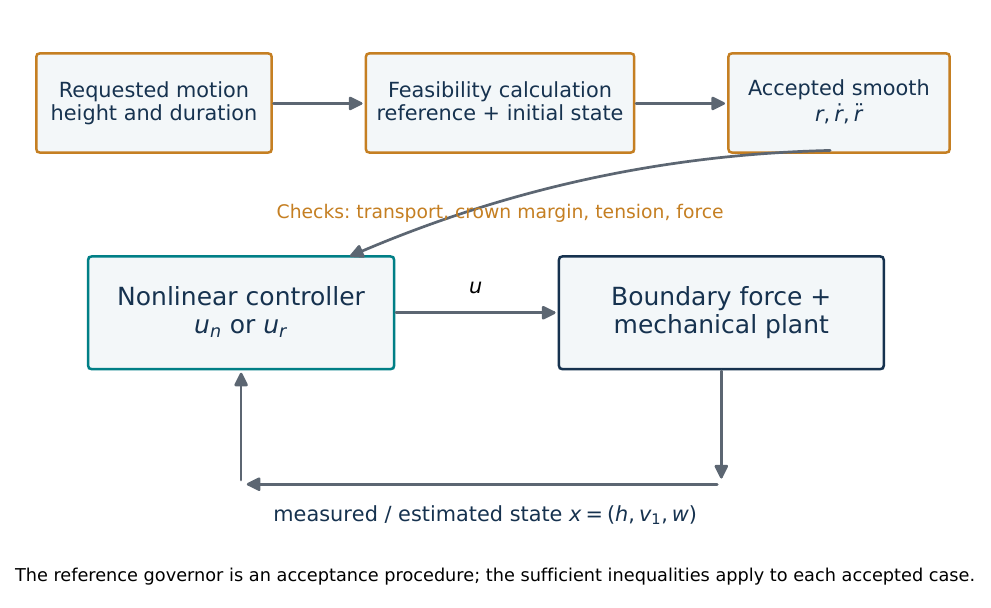}
\caption{The proposed control interpretation. The physical actuator applies the computed boundary force. Before a motion is requested, the reference is checked against the sufficient inequalities. The tracking results are proved for the unfiltered command. The optional crown filter is a separate, conditional construction in Section~\ref{sec:barrier}.}\label{fig:architecture}
\end{figure}

\subsection{What the mathematical programme is trying to establish}
The output is the height, but the plant has three states. Therefore an error estimate for $h-r$ is only one part of the argument. A useful result must also control the internal transport speed, keep the denominators away from zero, establish existence for the whole requested duration, and verify that the force corresponds to a tensile chain with an available actuator.

The logical target is the implication
\[
\begin{gathered}
\text{model assumptions + admissible initial data + explicit inequalities}\\
\Longrightarrow\ \text{a unique global model trajectory with tracking,}\\
\text{positive transport, crown margin, nonnegative tension and bounded input.}
\end{gathered}
\]
Each part of the conclusion requires its own argument. Error cancellation establishes tracking, while the falling-speed face inequalities control the internal state. Together, these estimates confine the trajectory to a compact set and allow the local solution to be continued. The force bounds then establish implementability within the ideal servo model. Under measurable bounded uncertainty, the same argument gives an explicit ultimate-error bound. With white noise the conclusion changes: one can still prove a stopped energy estimate, but a nonzero diffusion component normal to the crown boundary rules out almost-sure invariance from that boundary.

\subsection{Why this approach, and what is contributed to control theory}
A linear controller near the natural equilibrium is a reasonable engineering baseline, and the stable transfer function identifies its local sign and relative degree. Such a local description is insufficient for the present proof objective because positivity, tension and changing denominators impose nonlinear constraints. We use exact cancellation in the error coordinates to obtain a simple quadratic estimate while leaving one scalar internal equation whose boundary signs can be checked explicitly. This leads to sufficient conditions that can be evaluated directly and formalized in Lean.

The controller is not claimed to use force optimally, to stabilize the whole positive domain globally, or to be an experimentally validated design. Backstepping, comparison, invariant sets, control barrier functions and stochastic localization are established methods~\cite{krstic,khalil,ames,baudoin}. Here we connect these methods within the chain-fountain model to obtain a sufficient tracking region with an existence theorem and force checks, a quantitative uncertainty--force tradeoff, an exact compatibility condition for two particular robust constraints, and a sharp boundary obstruction for the actual mechanical diffusion. These results show how a performance calculation can be extended to a theorem that also establishes physical admissibility for a stable transport process.

For the phenomenon itself, the analysis separates the naturally selected height from the height maintained by boundary work and identifies limits that can be established before choosing the gains. It also distinguishes a robust description with bounded coefficients from a diffusion approximation. For control theory, the example shows why output stability, internal-state viability, actuator feasibility and solution construction require separate proofs. The formal development verifies the stated implications within the adopted model, whose experimental validation remains a separate question.

\subsection{Reading routes and the public mathematical companions}
Readers interested in mechanics can start with Sections~\ref{sec:model}--\ref{sec:spectrum}, then read the interpretation of the operating-region conditions in Section~\ref{sec:tube}. Readers interested in controller design can follow Sections~\ref{sec:tracking}--\ref{sec:barrier}. The white-noise branch in Section~\ref{sec:stochastic} is logically separate from the bounded-disturbance guarantee. The analytical tools used throughout are developed in Appendix~\ref{app:math}, and their connection to the checked declarations is recorded in Appendix~\ref{app:formal}.

Two public companions by the author provide further background. The MIACE companion~\cite{tomov13}, particularly Sections~2.5, 5.3, 5.6 and 7.1--7.2, covers absolute continuity, retained-state stability, quadratic energy and continuation. The 2DOF PIDF guide~\cite{tomov14}, Sections~2--4 and its reading map, explains how controller choices lead to explicit constraint geometry. Here feasibility is described by sufficient inequalities for the state, reference and force, rather than by a setpoint-weight polygon. The chain-fountain results and the stochastic transfer used here are proved in this manuscript and its own \href{https://github.com/luchesart/Chain-fountain/blob/main/Chain_Fountain_Lean_GitHub.zip}{formal companion}, not in either of those background documents.

\subsection*{Contributions and their scope}
The principal result concerns the combined tracking and admissibility problem. Theorem~\ref{thm:tube} and Proposition~\ref{prop:bounds} give explicit sufficient conditions on the initial errors, reference, gains and speed intervals. Under these conditions, exponential height tracking is accompanied by positive height and leg velocities, a crown-speed margin, nonnegative pickup tension and bounded force. Proposition~\ref{prop:global} constructs the unique global solution and proves continuation within this region. For a constant reference, Corollary~\ref{cor:fullstate} adds exponential convergence of the internal speed. The role of these results is to turn a calculation for the error coordinates into a statement about an admissible motion of the complete reduced plant.

Theorem~\ref{thm:robust} extends the construction to bounded measurable pickup uncertainty. Its error-energy envelope quantifies the benefit of additional damping, while Proposition~\ref{prop:bounds} quantifies the associated force demand. Thus disturbance rejection and mechanical feasibility can be assessed together. The damping method is established, but its consequences here depend on the fountain's pickup law and internal-speed bounds.

Several supporting results explain where this construction applies. The Routh--Hurwitz factorization and the local Lyapunov argument in Theorem~\ref{thm:spectrum} establish nonlinear local exponential stability over the stated parameter range and the transfer calculation identifies the stable real zero. Propositions~\ref{prop:energy} and~\ref{prop:setpoint} identify actuator power and the steady heights accessible through braking or assistance. Proposition~\ref{prop:barrier} gives an exact compatibility condition for the robust crown and pickup inequalities. It is a test for those two constraints, with actuator limits imposed by interval intersection. The counterexamples explain why it does not characterize the whole viable domain.

The stochastic results address a different disturbance model. Section~\ref{sec:stochastic} constructs the maximal local mechanical solution, proves a stopped energy estimate and transfers a small-time integral estimate to its actual crown coordinate. The resulting sharp law of the iterated logarithm rules out almost-sure crown invariance from a boundary state with nonzero normal diffusion. This provides a model-specific pathwise account of an obstruction already familiar from stochastic viability theory.

Section~\ref{sec:zerog} treats open transport and closed circulation separately. In the open zero-gravity model it proves a finite-horizon obstruction to uniformly positive flow and exhibits an exact decaying-flow solution. For a closed string, it starts from a classical fixed-shape circulation family and constructs speed tracking under a distributed tangential force, including measurable disturbances and physical margins. The conclusion concerns circulation within that family. Spatial shape stability requires additional equations and estimates.

The numerical examples illustrate the reduced models and make the force and motion bounds inspectable. The Lean companion checks the stated implications and records their dependencies. These forms of evidence complement the analytical argument, while testing the adopted mechanics in an actuated chain requires independent experiments or discrete-chain calculations.

\section{Model, input, and physical domain}\label{sec:model}
Write $H=h_1>0$ for the drop from the pickup level to the floor, $h=h_2>0$ for the apex height above pickup, $v_1>0$ for the upward speed, and $w=v_2>0$ for the downward speed. The chain has linear density $\lambda>0$. Both legs are vertical and carry uniform axial velocity. The crown has negligible mass and length in comparison with the two legs. Pickup and deposition are represented by instantaneous constitutive laws. Elastic waves, finite crown radius, lateral motion, and changing pile elevation are omitted.

Conservation of airborne chain length gives
\begin{equation}
 \dot h=\frac{v_1-w}{2}=:q,\qquad
 s:=\frac{v_1+w}{2}=w+q,\qquad T_C=\lambda s^2.\label{eq:kin}
\end{equation}
The different leg velocities are consistent with inextensibility: relative to the moving crown, both legs pass through the turn at speed $s$.

\subsection{An ideal boundary force}
Let an actuator exert a downward force $\lambda u$ in the pickup zone. Thus $u>0$ denotes braking and $u<0$ denotes upward assistance. If the upward pile reaction is $\alpha\lambda v_1^2$, the pickup momentum balance is
\begin{equation}
 T_P+\alpha\lambda v_1^2-\lambda u=\lambda v_1^2,
 \qquad T_P=\lambda\bigl[(1-\alpha)v_1^2+u\bigr].\label{eq:pickup}
\end{equation}
A driven traction device is one possible realization of a signed force. A passive brake can provide only $u\ge0$. An exit-angle change is not automatically equivalent to this input, because it also changes the geometry and possibly $\alpha$. Equation~\eqref{eq:pickup} assumes that the local force acts on the transported chain without otherwise changing the effective pickup law. This assumption requires calibration in a physical apparatus.

At the lower boundary, $T_F=\lambda\beta w^2$ denotes the tensile force pulling downward on the falling leg. With $S=s^2$, the controlled model is
\begin{subequations}\label{eq:model}
\begin{align}
 \dot h&=(v_1-w)/2,\\
 \dot v_1&=\frac{S-(1-\alpha)v_1^2-gh-u}{h}=:f_1(x,u,\alpha),\\
 \dot w&=g+\frac{-S+\beta w^2}{H+h}=:f_2(x),
\end{align}
\end{subequations}
where $x=(h,v_1,w)$ and
\begin{equation}
 \D=\{(h,v_1,w):h>0,\ v_1>0,\ w>0\}.\label{eq:domain}
\end{equation}
Forward pickup is part of the domain: continuing the same pickup law after $v_1$ changes sign has no stated physical justification.

\begin{assumption}\label{ass:parameters}
The nominal coefficients satisfy $0<\alpha_0<1/2$, $\beta\ge0$, and $\alpha_0+\beta<1/2$. In the uncertain model, $\alpha(t)=\alpha_0+d(t)$, where $d$ is measurable and $|d(t)|\le\bar d$, with $0<\alpha_-:=\alpha_0-\bar d$ and $\alpha_+:=\alpha_0+\bar d<1/2-\beta$.
\end{assumption}
We adopt this parameter range from the reduced-model setting of \cite{bigginswarner,biggins}, without assuming that it covers the full possible range for every chain. The input has units $\mathrm{m^2\,s^{-2}}$ and is constrained by
\begin{equation}
 u_{\mathrm{act},-}\le u\le u_{\mathrm{act},+}.\label{eq:actlimits}
\end{equation}
The theorems first concern an ideal force servo. A force lag changes the state equations and is examined only as a numerical sensitivity test in Section~\ref{sec:numerics}.

For a finite chain, the remaining pickup supply $\ell$ obeys $\dot\ell=-v_1$. Physical conclusions apply only while $\ell>0$ and the deposition and geometry assumptions hold. Statements for all $t\ge0$ concern the idealized model with sustained supply. A finite experimental chain can follow such a motion only until its supply is exhausted.

\subsection{Recovering the equations from transport balances}\label{sec:derive}
The total airborne length in this idealization is $L=H+2h$. The pickup boundary introduces chain at rate $v_1$ and the deposition boundary removes it at rate $w$. Thus $\dot L=v_1-w$, which gives the first equation in~\eqref{eq:kin}. At a crown moving upward at speed $q$, the rising chain arrives at relative speed $v_1-q$, while the descending chain leaves at relative speed $w+q$. Both equal $s=(v_1+w)/2$. Turning this flux through an ideal massless crown gives the effective tension $T_C=\lambda s^2$ of the adopted model.

\begin{figure}[htbp]
\centering\includegraphics[width=.96\linewidth]{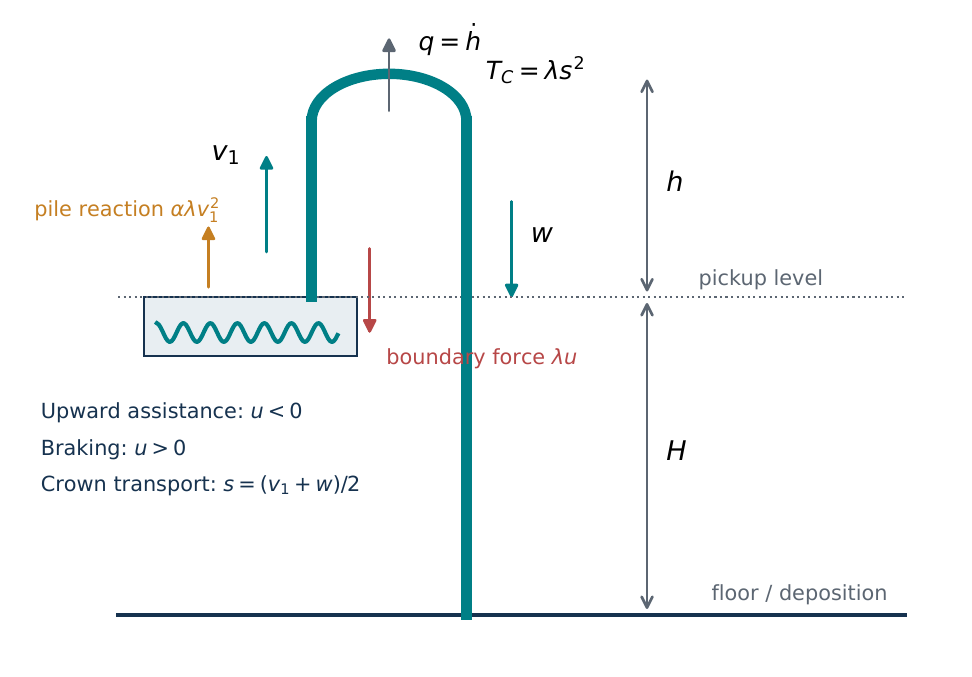}
\caption{Coordinates and force signs. The two legs are separated for legibility in this schematic drawing, which is not to scale. The analytical crown has negligible length and mass. The proposed traction force is an additional boundary action, distinct from the effective upward pile reaction.}\label{fig:geometry}
\end{figure}

To make the moving-boundary terms explicit, we write the momentum balance for each leg. Applying Newton's law as though $\lambda h$ were a fixed mass would omit these terms. With upward momentum taken as positive, the rising leg has momentum $\lambda h v_1$. Its incoming momentum flux at pickup is $\lambda v_1^2$, and the outgoing flux at the moving crown is $\lambda s v_1$. Hence
\begin{align*}
 \frac{\dd}{\dd t}(\lambda h v_1)
 &=T_C-T_P-\lambda gh+\lambda v_1^2-\lambda s v_1\\
 &=T_C-T_P-\lambda gh+\lambda qv_1.
\end{align*}
The left side is $\lambda h\dot v_1+\lambda qv_1$, so the transport terms cancel and
\[
\lambda h\dot v_1=T_C-T_P-\lambda gh.
\]
Substituting $T_P=\lambda[(1-\alpha)v_1^2+u]$ gives the rising-leg equation in~\eqref{eq:model}.

The pickup law is a constitutive assumption in this calculation. Pantaleone's measurements of fountain height and speed, together with a model of how a link lifts off from the container, show how that boundary response can be examined experimentally~\cite{pantaleone2017}. Such a comparison is needed to assign a physical value to $\alpha$ for a particular chain and packing. The momentum balance alone does not determine it.

For the falling leg, take downward momentum as positive. The leg has length $H+h$ and momentum $\lambda(H+h)w$. The crown supplies momentum flux $\lambda sw$, and the floor removes $\lambda w^2$, so
\begin{align*}
\frac{\dd}{\dd t}\{\lambda(H+h)w\}
 &=\lambda g(H+h)+T_F-T_C+\lambda sw-\lambda w^2\\
 &=\lambda g(H+h)+T_F-T_C+\lambda qw.
\end{align*}
Again the changing-length term cancels, giving
\[
\lambda(H+h)\dot w=\lambda g(H+h)+T_F-T_C.
\]
The choice $T_F=\lambda\beta w^2$ yields the final equation. A stationary floor can exert an upward contact force while the chain immediately above it experiences a downward tensile pull. The arriving links are brought to rest partly through their interaction with the still-moving chain. Experiments and models of falling chains demonstrate how this interaction can accelerate the remaining chain beyond free fall~\cite{hamm2010,grewal2011}. These results motivate allowing a positive $\beta$, although its value and dependence on impact conditions must be determined for the apparatus being modeled.

With $u=0$, the three equations are Biggins's transient vertical-fountain equations~\cite[Eqs.~(10)--(12)]{biggins}. The calculation above makes their transport terms explicit and adds the prescribed boundary force. It recovers the reduced equations from the assumed boundary laws. Deriving the constitutive coefficients from link-scale mechanics is outside its scope.

\subsection{States, units, and what an implementation would measure}
\begin{table}[htbp]\centering\small
\begin{tabular}{>{\raggedright\arraybackslash}p{.16\linewidth}p{.20\linewidth}p{.53\linewidth}}\toprule
Quantity & SI unit & Meaning in the argument\\\midrule
$h,H,r$ & m & apex height, fixed drop, desired height\\
$v_1,w,s,q$ & m\,s$^{-1}$ & rising, falling, crown-relative and apex speeds\\
$\lambda$ & kg\,m$^{-1}$ & chain mass per unit length\\
$\alpha,\beta,d$ & dimensionless & effective boundary coefficients and perturbation\\
$u$ & m$^2$\,s$^{-2}$ & downward boundary force divided by $\lambda$\\
$E$ & J & mechanical storage, after a constant shift\\
$V$ & m$^2$\,s$^{-2}$ & weighted error energy, not mechanical energy\\
$k_1,k_2;k_0$ & s$^{-1}$; s$^{-2}$ & error-dynamics gains\\
$\kappa_d;\sigma$ & s$^3$\,m$^{-2}$; s$^{1/2}$ & robust gain; diffusion intensity\\\bottomrule
\end{tabular}
\caption{Dimensions expose the distinction between physical storage and a control Lyapunov function. Brownian motion is assigned units $\mathrm{s}^{1/2}$ in the dimensional SDE.}\label{tab:units}
\end{table}

The ideal feedback requires the height and both leg speeds, which could in principle be obtained from imaging and speed measurements. No observer or sampling theorem is part of the present result. The command also uses the nominal coefficients and a prescribed $C^2$ reference, whose derivatives are known functions rather than estimates from noisy measurements. The formula for $\dot a$ later uses the model to avoid differentiating a measured velocity.

Actuator lag, measurement error, a changing $H$, and depletion of the pile change either the equations or the admissible duration of the motion. The same applies if a command is saturated after it has been computed without first establishing its feasibility. We therefore check the force inequalities before accepting a reference. The closed-loop theorem assumes an ideal servo, and any departure from that assumption in the proposed apparatus requires further analysis.

\FloatBarrier\section{Energy and achievable steady heights}\label{sec:energy}
The kinetic energy of the airborne chain is $\lambda[h v_1^2+(H+h)w^2]/2$. Its gravitational potential energy, relative to pickup level, is $\lambda g(h^2-H^2/2)$. Since $H$ is fixed, a constant shift yields the nonnegative storage
\begin{equation}
 E(x)=\frac{\lambda}{2}\bigl[h v_1^2+(H+h)w^2\bigr]+\lambda gh^2.\label{eq:storage}
\end{equation}
\begin{proposition}[Energy balance within the reduced model]\label{prop:energy}
Along \eqref{eq:model}, for constant or measurable bounded $\alpha$,
\begin{equation}
 \dot E=\lambda gHw-\lambda(\tfrac12-\alpha)v_1^3
 -\lambda(\tfrac12-\beta)w^3-\lambda u v_1.\label{eq:energy}
\end{equation}
For measurable coefficients, this identity holds almost everywhere.
\end{proposition}
\begin{proof}
Differentiating \eqref{eq:storage} and substituting the leg equations gives, after division by $\lambda$,
\begin{align*}
 \dot E/\lambda={}&\tfrac12\dot h(v_1^2+w^2)
 +v_1[S-(1-\alpha)v_1^2-gh-u]\\
 &+w[g(H+h)-S+\beta w^2]+2gh\dot h.
\end{align*}
The gravitational terms reduce to $gHw$. Using $\dot h=(v_1-w)/2$ and $S=(v_1+w)^2/4$, the remaining crown and changing-mass terms are
\[
 (v_1-w)S+\tfrac14(v_1-w)(v_1^2+w^2)
 =\tfrac12(v_1^3-w^3).
\]
Collecting the cubic and input terms proves the identity.
\end{proof}
The supply $\lambda gHw$ comes from transport to the lower level. The next two terms are nonnegative losses under Assumption~\ref{ass:parameters}, and $\lambda u v_1$ is power extracted by the actuator when $u>0$. The gravitational supply is essential to this energy balance, so the identity alone does not establish passivity with respect to $u$.

\begin{proposition}[Steady-state input map]\label{prop:setpoint}
For a constant coefficient $\alpha$ and prescribed height $r>0$, the positive steady state requires
\begin{equation}
 v_r^2=\frac{g(H+r)}{1-\beta},\qquad
 u_r=\frac{g[\alpha H-(1-\alpha-\beta)r]}{1-\beta}.\label{eq:setpoint}
\end{equation}
Its pickup tension is strictly positive:
\begin{equation}
 T_{P,r}/\lambda=\frac{g(H+\beta r)}{1-\beta}>0.\label{eq:steady_tension}
\end{equation}
Writing
\begin{equation}
 R(u)=\frac{\alpha H-(1-\beta)u/g}{1-\alpha-\beta},
\end{equation}
the heights feasible under \eqref{eq:actlimits} are precisely
\begin{equation}
 (0,\infty)\cap[R(u_{\mathrm{act},+}),R(u_{\mathrm{act},-})]
\end{equation}
within this steady reduced model.
\end{proposition}
\begin{proof}
At equilibrium $v_1=w=v_r$. The falling-leg equation gives the first identity, and the rising-leg equation gives $u_r=\alpha v_r^2-gr$. Substitution into \eqref{eq:pickup} gives \eqref{eq:steady_tension}. The map $u_r$ is strictly decreasing in $r$, so the actuator interval is inverted in the stated order.
\end{proof}
In particular, with an unrestricted passive brake $u\ge0$, positive steady heights satisfy
\begin{equation}
 0<r\le h^*,\qquad h^*:=\frac{\alpha H}{1-\alpha-\beta}.\label{eq:brakeheight}
\end{equation}
Raising the sustained fountain above $h^*$ requires active assistance. This limitation follows from the force balance before a controller is selected. For an uncertain coefficient and a known desired height, the monotonic increase of $u_r(\alpha)$ with $\alpha$ gives a necessary robust actuator test: its values at both coefficient endpoints must lie within the actuator interval.

\subsection{What the energy and equilibrium results tell us}
As the apex moves, the masses of the two legs change, producing the changing-mass terms in the energy derivative. The boundary powers are cubic because a force proportional to speed squared acts on chain moving at that speed. The constant shift in $E$ changes neither its derivative nor the force balance. The error energy $V$ introduced later serves a different purpose: it measures distance from a requested motion and is designed to decay.

The loss terms in~\eqref{eq:energy} describe the net energy transfer out of the airborne chain under the adopted pickup and deposition laws. Their derivation must account for the energy carried by material crossing each boundary as well as for the work of the tensions. Virga's treatment of the fountain through dissipative shocks provides a more general framework for imposing momentum and energy balances at such discontinuities~\cite{virga2014}. Here the two effective coefficients specify the boundary behavior before the energy identity is derived. The identity therefore tests consistency of the reduced model without establishing a microscopic law for the losses.
\begin{figure}[htbp]\centering
\includegraphics[width=\linewidth]{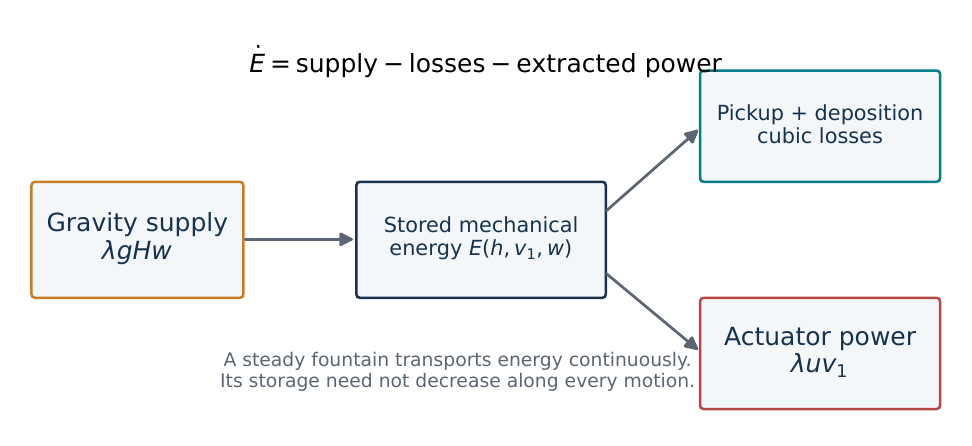}
\caption{Mechanical energy accounting. The actuator term is extraction when $u>0$ and injection when $u<0$. A stationary fountain is an open transport process with nonzero power throughput.}\label{fig:energyflow}
\end{figure}

At a steady state the storage derivative is zero, although gravity continues to supply power and the two boundary mechanisms continue to dissipate it. The condition $E\ge0$ alone therefore cannot establish stability of the unforced fountain, since the right side of~\eqref{eq:energy} includes a positive supply. Local stability requires a separate proof. The energy identity does, however, establish the physical meaning of the input sign before we use it in a stability argument.

At a requested constant height, the falling equation fixes $v_r$ independently of the pickup coefficient, after which the rising equation determines $u_r$. Height and transport speed therefore cannot be selected independently. The steady command has sensitivities
\[
\frac{\dd u_r}{\dd r}=-\frac{g(1-\alpha-\beta)}{1-\beta}<0,
\qquad
\frac{\partial u_r}{\partial\alpha}=\frac{g(H+r)}{1-\beta}=v_r^2>0.
\]
Increasing the desired height requires less braking or more assistance. Increasing the pickup reaction at the same height requires more braking. These derivative signs also explain the reversed endpoints of the feasible height interval and why two coefficient endpoints suffice for the uncertain steady-input test.

\begin{figure}[htbp]\centering
\includegraphics[width=\linewidth]{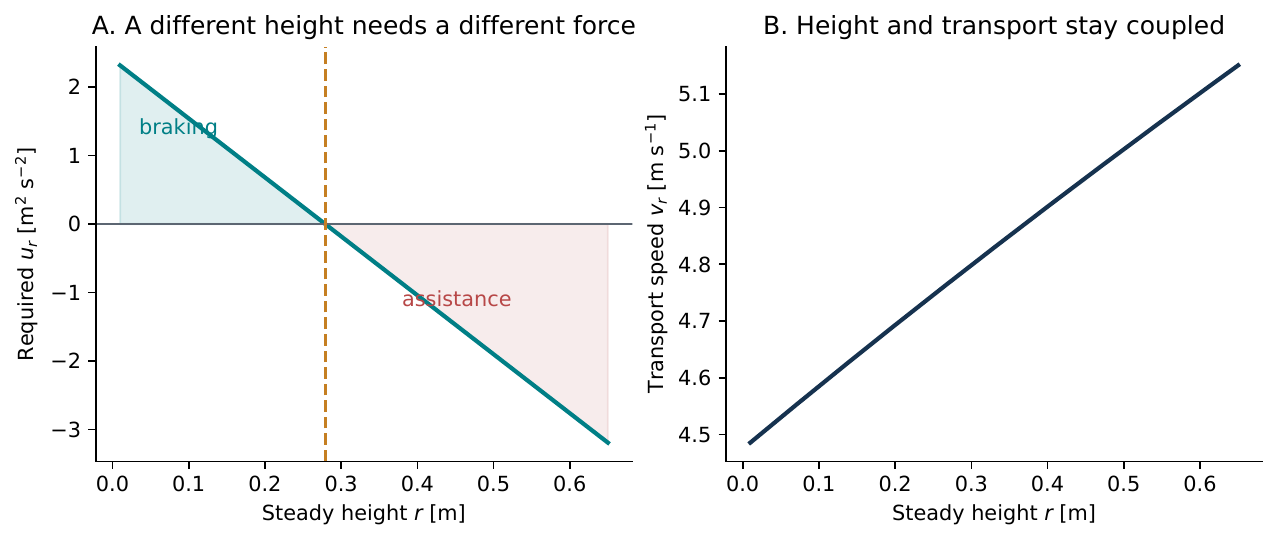}
\caption{Analytical equilibrium maps for $H=2$, $g=9.81$, $\alpha=0.12$, $\beta=0.02$. The dashed line is the natural height $h^*\simeq0.27907$ m. These curves are exact evaluations of~\eqref{eq:setpoint}, not simulation data.}\label{fig:steadymap}
\end{figure}

For example, the same parameters give $u_r\simeq0.68069\,\mathrm{m^2s^{-2}}$ at $r=0.20$ m and $u_r\simeq-0.61062\,\mathrm{m^2s^{-2}}$ at $r=0.35$ m. Natural stability at $0.27907$ m does not remove the steady force needed to maintain either of these heights. Multiplying by the measured linear density $\lambda$ gives the physical force. We do not fit a value of $\lambda$ here.

\FloatBarrier\section{Local stability and minimum-phase structure}\label{sec:spectrum}
For $u=0$ and $\alpha=\alpha_0$, the unique positive equilibrium is
\begin{equation}
 x^*=(h^*,v^*,v^*),\quad h^*=\frac{\alpha_0H}{1-\alpha_0-\beta},\quad
 v^*=\sqrt{\frac{gH}{1-\alpha_0-\beta}}.
\end{equation}
Set $\omega=v^*/h^*$ and $\rho=h^*/(H+h^*)=\alpha_0/(1-\beta)$.
\begin{theorem}\label{thm:spectrum}
Under Assumption~\ref{ass:parameters} with $d=0$, $x^*$ is locally exponentially stable. The linearized transfer function from $\Delta u$ to $\Delta h$ is
\begin{equation}
 W(p)=-\frac{1}{2h^*}\frac{p+2\alpha_0\omega}{p^3+a_2p^2+a_1p+a_0},\label{eq:transfer}
\end{equation}
where
\begin{align}
 a_2&=(1-\rho)\omega,\nonumber\\
 a_1&=\left[\tfrac52\alpha_0(1+\rho)-4\alpha_0^2\right]\omega^2,\label{eq:coeff}\\
 a_0&=\alpha_0^2(1-\rho)\omega^3.\nonumber
\end{align}
The transfer function has relative degree two, negative static gain, and one finite zero $-2\alpha_0\omega=-2g/v^*<0$.
\end{theorem}
\begin{proof}
The Jacobian and input vector are
\begin{equation}
 J=\begin{bmatrix}
 0&1/2&-1/2\\
 -\alpha_0\omega^2&-(1-2\alpha_0)\omega&\omega\\
 \alpha_0\rho\omega^2&-\rho\omega&-(1-2\beta)\rho\omega
 \end{bmatrix},\qquad
 B_u=\begin{bmatrix}0\\-1/h^*\\0\end{bmatrix}.\label{eq:jac}
\end{equation}
Expanding the characteristic polynomial and using $\alpha_0=(1-\beta)\rho$ gives \eqref{eq:coeff}. Each coefficient is positive, and the remaining cubic Routh--Hurwitz inequality follows from the factorization
\begin{equation}
 a_2a_1-a_0=\tfrac52\alpha_0(1-\rho)(1+\rho-2\alpha_0)\omega^3>0.\label{eq:routh}
\end{equation}
Indeed $0<\rho<1$ and $1+\rho-2\alpha_0>0$. To pass from this spectral calculation to nonlinear stability, the Lean companion constructs an explicit positive quadratic form for the cubic companion matrix and transfers it to the actual Jacobian through an injective observability map. A bound on the nonlinear derivative remainder then gives a neighbourhood in which the quadratic form decreases. Using a bounded extension and a first-exit argument, the formal proof constructs a global classical solution for every sufficiently close initial state, with $\|x(t)-x^*\|\le C e^{-ct}\|x(0)-x^*\|$ and positive height for all $t\ge0$. It also proves uniqueness against every classical competitor on any finite forward interval, without requiring that competitor to have positive height or remain in the neighbourhood. Finally, with $C=(1,0,0)$, direct evaluation of $C(pI-J)^{-1}B_u$ gives \eqref{eq:transfer}. Its numerator has the stated sign and root. Substitution of $p=-2\alpha_0\omega$ into the denominator gives $-10\alpha_0^2\rho\omega^3\ne0$, which rules out pole-zero cancellation.
\end{proof}
The passage from a Hurwitz linearization to a quadratic Lyapunov estimate is standard nonlinear stability theory~\cite[Chapter 4]{khalil}. Here the parameter factorization and the construction for the mechanical Jacobian supply its hypotheses explicitly.

For $H=2\,\mathrm m$, $\alpha_0=0.12$, $\beta=0.02$ and $g=9.81\,\mathrm{m\,s^{-2}}$, the poles are $-0.92112683$ and $-7.04926243\pm4.36925499i$, in $\mathrm{s^{-1}}$. The zero is $-4.10770009\,\mathrm{s^{-1}}$, and \eqref{eq:routh} is approximately $1164.78479\,\mathrm{s^{-3}}$. This is a local stability result for the three-state equilibrium. Global attraction and transverse-chain modes require separate analysis, with the latter lying outside the model.

\subsection{Reading the stability calculation mechanically}
The first row of $J$ expresses the kinematic relation between apex motion and the difference in leg speeds. The entry $J_{21}= -\alpha_0\omega^2$ is obtained by differentiating the rising acceleration with respect to height. At equilibrium, the numerator of that acceleration vanishes, so the derivative comes from the height dependence of the force balance. The entry $J_{31}=\alpha_0\rho\omega^2$ has the opposite sign. The speed couplings $J_{23}=\omega$ and $J_{32}=-\rho\omega$ arise from the crown tension shared by the two legs. Stability depends on how these terms act together.

The input derivative is $B_u=(0,-1/h^*,0)^\top$. A larger downward force first reduces rising acceleration and then changes the height through the kinematic relation. This gives the input/output model relative degree two. Its negative DC gain agrees with the equilibrium derivative:
\[
W(0)=-\frac{1-\beta}{g(1-\alpha_0-\beta)}=\frac{\dd R}{\dd u}.
\]
The zero is strictly in the left half-plane, so this local input/output model has no unstable finite zero. This supports the local tracking design. To determine whether a finite excursion is physically admissible, we must also examine the internal nonlinear state, which is identified in Section~\ref{sec:tracking}.

The factored Hurwitz determinant establishes stability throughout the stated parameter region. For a cubic, positivity of the coefficients must be accompanied by the inequality $a_2a_1>a_0$. This algebraic test places the eigenvalues of $J$ in the open left half-plane. The passage to nonlinear local stability then uses a positive definite quadratic Lyapunov function for the displacement from equilibrium. Its derivative is negative for the linearized system, and remains negative in a sufficiently small neighbourhood after the nonlinear remainder is included. Appendix~\ref{app:hurwitz} explains this argument and its formal realization. The resulting local neighbourhood and the invariant tracking region constructed below serve different purposes. The former describes the natural stability of the equilibrium, while the latter supports controlled motion under explicit physical constraints.

\section{Dimensionally consistent tracking and internal dynamics}\label{sec:tracking}
We use the backstepping construction of \cite{krstic} and choose its weights so that the tracking-error measure is dimensionally consistent. Let $r\in C^2([0,\infty))$ satisfy
\begin{equation}
 0<r_-\le r(t)\le r_+,\qquad |\dot r|\le R_1,\qquad |\ddot r|\le R_2.\label{eq:refbounds}
\end{equation}
Choose $k_1,k_2>0$ with units $\mathrm{s^{-1}}$ and $k_0>0$ with units $\mathrm{s^{-2}}$. We define the height error $e$, the virtual velocity $a$ and its mismatch $z$, and introduce the positive definite quadratic Lyapunov function $V$ for the tracking-error system:
\begin{equation}
 e=h-r,\quad a=w+2\dot r-2k_1e,\quad z=v_1-a,\quad
 V=\frac{k_0}{2}e^2+\frac18z^2.\label{eq:errors}
\end{equation}
Both terms in $V$ have units $\mathrm{m^2\,s^{-2}}$. In particular, $2V$ is the squared Euclidean norm of the weighted error vector $(\sqrt{k_0}e,z/2)$, whose components both have velocity units. Thus $V$ measures the size of the two tracking errors on a common scale and vanishes exactly when $e=z=0$. The virtual velocity derivative is available from the model:
\begin{equation}
 \dot a=f_2(x)+2\ddot r+2k_1^2e-k_1z.\label{eq:adot}
\end{equation}
This avoids numerical differentiation of measured velocities, although it still assumes state measurements or estimates and a sufficiently accurate model.

\begin{proposition}[Nominal error cancellation]\label{prop:tracking}
For $\alpha=\alpha_0$, the control
\begin{equation}
 \un=S-(1-\alpha_0)v_1^2-gh+h(2k_0e-\dot a+k_2z)\label{eq:nominal}
\end{equation}
gives, on every interval of existence in $\D$ on which this input is applied,
\begin{equation}
 \dot e=-k_1e+\tfrac12z,\qquad \dot z=-2k_0e-k_2z,\qquad
 \dot V=-k_0k_1e^2-\tfrac14k_2z^2\le-2\kappa V,\label{eq:errdyn}
\end{equation}
where $\kappa=\min(k_1,k_2)$. Thus $V(t)\le V(0)e^{-2\kappa t}$ on that interval.
\end{proposition}
\begin{proof}
The first equation follows from \eqref{eq:kin} and \eqref{eq:errors}. Substituting \eqref{eq:nominal} into the rising-leg equation and subtracting \eqref{eq:adot} gives the second. In $\dot V=k_0e\dot e+z\dot z/4$, the mixed terms cancel, yielding the stated estimate.
\end{proof}
The controller is equivalently an exact output-acceleration assignment:
\begin{equation}
 \dot q=\ddot r-(k_1+k_2)(q-\dot r)-(k_0+k_1k_2)e=:A_0.\label{eq:acceleration}
\end{equation}
This form shows the dimensional meaning of the gains and how the reference acceleration enters the command. Whether a chosen reference and initial state are physically feasible still depends on the operating conditions derived below.

The remaining state obeys
\begin{equation}
 \dot w=\frac{g(H+h)-(1-\beta)w^2-2qw-q^2}{H+h}.\label{eq:internal}
\end{equation}
For a constant reference on the exact tracking manifold, $h=r$ and $q=0$. Then
\begin{equation}
 \dot w=-\frac{1-\beta}{H+r}(w+w_r)(w-w_r),\qquad
 w_r=\sqrt{\frac{g(H+r)}{1-\beta}}.\label{eq:zero}
\end{equation}
The positive half-line is invariant and $|w(t)-w_r|\le|w(0)-w_r|e^{-ct}$ with $c=(1-\beta)w_r/(H+r)>0$. This proves global exponential attraction on the positive branch of the constant-reference zero dynamics. For a varying $q$, however, the boundary derivative at $w=0$ is $\dot w=g-q^2/(H+h)$, which can be negative even when $q$ is bounded. Unrestricted ISS of \eqref{eq:internal} with respect to $q$ therefore does not follow. The operating region below addresses this difficulty by constraining the signals jointly.

\subsection{Designing the error coordinates, one cancellation at a time}
The quadratic Lyapunov function $V$ in~\eqref{eq:errors} guides the tracking design. Since $k_0>0$, it is positive for every nonzero error pair $(e,z)$ and vanishes at exact tracking. We choose the feedback so that its derivative measures decay of both errors.

The kinematics give the first step. To obtain $\dot e=-k_1e$, the rising speed would have to equal $a=w+2\dot r-2k_1e$. The difference $z=v_1-a$ between the actual and prescribed rising speeds then gives
\[
\dot e=(v_1-w)/2-\dot r=-k_1e+z/2.
\]
Differentiating the prescribed virtual speed gives
\[
\dot a=\dot w+2\ddot r-2k_1\dot e
=f_2+2\ddot r+2k_1^2e-k_1z.
\]
To obtain the second error equation, we prescribe $\dot v_1=\dot a-2k_0e-k_2z$ and solve the rising force balance for $u$. This gives~\eqref{eq:nominal}. The coefficient $2k_0$ makes the two mixed products in $\dot V$ cancel:
\[
\underbrace{k_0e(z/2)}_{\text{from height}}+
\underbrace{(z/4)(-2k_0e)}_{\text{from velocity}}=0.
\]
The remaining terms give $\dot V=-k_0k_1e^2-k_2z^2/4$, which is negative whenever $(e,z)\ne(0,0)$. More precisely,~\eqref{eq:errdyn} gives $\dot V\le-2\kappa V$. Integrating this inequality proves exponential decay on every interval where the feedback is applied and the solution stays in the physical domain. In weighted coordinates $\xi=(\sqrt{k_0}e,z/2)$, we have $V=\|\xi\|^2/2$, so the estimate directly bounds the tracking errors in a dimensionally consistent norm. The weights in $V$ thus serve two purposes: they put the errors on a common physical scale and match the coupling terms in the error equations.

Eliminating $z=2(\dot e+k_1e)$ gives the familiar second-order equation
\[
\ddot e+(k_1+k_2)\dot e+(k_0+k_1k_2)e=0.
\]
The gains determine the damping and stiffness coefficients explicitly. The resulting convergence estimate allows overshoot and nonmonotonic height motion for some initial errors. The monotone-tracking results of the MIACE companions concern a different controller and require their own hypotheses.

\subsection{Why the internal equation cannot be skipped}
The coordinate change is invertible at each time:
\[
h=r+e,\qquad q=\dot r-k_1e+z/2,\qquad v_1=w+2q.
\]
In these coordinates the nominal system is triangular. The two errors evolve autonomously, while $w$ is driven by $h$ and $q$. The Lyapunov estimate controls $(e,z)$, but leaves the internal speed $w$ to be bounded separately. On exact constant tracking, its scalar equation factors as~\eqref{eq:zero}. Away from that manifold, the forcing terms $-2qw-q^2$ can drive $w$ toward zero even as the output error decreases.

The invariant region constrains $q$ sufficiently to make the falling-speed vector field point inward at two chosen faces. The inverse coordinate relation then also bounds $v_1$. Combining these speed bounds with the Lyapunov estimate gives the physical admissibility result: the output tracks the reference while the internal transport remains within the stated operating limits.

\section{An invariant physical tracking region}\label{sec:tube}
Fix a proposed error-energy bound $\Vbar>0$ and define
\begin{align}
 c_e&=\sqrt{2\Vbar/k_0},& c_z&=\sqrt{8\Vbar},&
 c_q&=\sqrt{2\Vbar}\sqrt{1+k_1^2/k_0},\label{eq:cdefs}\\
 h_-&=r_--c_e,&h_+&=r_++c_e,&Q&=R_1+c_q.\nonumber
\end{align}
On $V\le\Vbar$, Cauchy--Schwarz gives $|e|\le c_e$, $|z|\le c_z$, and $|q-\dot r|\le c_q$, so $h\in[h_-,h_+]$ and $|q|\le Q$. Choose falling-speed bounds $0<w_-<w_+$.

\begin{theorem}[Nominal invariant tracking region]\label{thm:tube}
Suppose $h_->0$, the initial state satisfies $V(0)\le\Vbar$ and $w(0)\in[w_-,w_+]$, and
\begin{align}
 w_-&>2Q, &w_--Q&\ge v_{\min}>0,\label{eq:positive}\\
 g(H+h_-)&\ge(1-\beta)w_-^2+2Qw_-+Q^2,\label{eq:lowerface}\\
 g(H+h_+)&\le(1-\beta)w_+^2-2Qw_++Q^2.\label{eq:upperface}
\end{align}
Under \eqref{eq:nominal}, with constant $\alpha=\alpha_0$, the ideal model has a solution for all $t\ge0$ and
\begin{equation}
 V(t)\le\Vbar,\quad h_-\le h(t)\le h_+,\quad
 w_-\le w(t)\le w_+,\quad v_1(t)>0,\quad s(t)\ge v_{\min}.\label{eq:tube}
\end{equation}
The error estimate in Proposition~\ref{prop:tracking} holds for all $t$. Physical implementation additionally requires pickup-tension and actuator feasibility, supplied by Proposition~\ref{prop:bounds} below.
\end{theorem}
\begin{proof}
Consider a solution up to its first possible exit from the proposed region. By \eqref{eq:errdyn}, $V$ cannot increase through $\Vbar$. The bounds in \eqref{eq:cdefs} therefore apply. Since $w_->2Q$, the function $(w+q)^2$ is increasing in $q\in[-Q,Q]$ at either falling-speed face. At $w=w_-$ the numerator of \eqref{eq:internal} is at least the left side minus the right side of \eqref{eq:lowerface}, hence is nonnegative. At $w=w_+$ it is at most the left side minus the right side of \eqref{eq:upperface}, hence is nonpositive. The speed interval is invariant. Moreover $v_1=w+2q\ge w_--2Q>0$ and $s=w+q\ge w_--Q\ge v_{\min}$. These estimates hold on each finite prefix of the local solution. The positive lower bound on $h$ and the velocity bounds place that solution in a fixed compact subset of the regularity domain. Proposition~\ref{prop:global} supplies the local construction, uniqueness, and continuation argument, using only these finite-prefix estimates. The tracking guarantee therefore includes construction of the global solution from the initial data.
\end{proof}

\begin{corollary}\label{cor:fullstate}
Under Theorem~\ref{thm:tube}, if $r$ is constant, the full state converges exponentially to $(r,w_r,w_r)$.
\end{corollary}
\begin{proof}
Let $\eta=w-w_r$. Subtracting the constant-reference balance gives
\[
 \dot\eta=-\frac{(1-\beta)(w+w_r)}{H+h}\eta
 +\frac{ge-2qw-q^2}{H+h}.
\]
The decay coefficient is bounded below by $(1-\beta)(w_-+w_r)/(H+h_+)>0$. The remaining terms decay exponentially because $e,z$, and hence $q$, do so and $w$ is bounded. Variation of constants gives exponential decay at any positive rate smaller than the minimum of the two decay rates. Finally $v_1=w+2q$ has the same limit.
\end{proof}

\subsection{Explicit tension and actuator bounds}
The invariant region concerns the state equations. Its physical use also needs $T_P\ge0$ and a force demand within \eqref{eq:actlimits}. The following conservative bounds make those checks finite algebraic calculations.

Set
\begin{align}
 v_-&=w_--2Q,&v_+&=w_++2Q,&s_-&=w_--Q,&s_+&=w_++Q,\nonumber\\
 A&=R_2+(k_1+k_2)c_q+(k_0+k_1k_2)c_e.\label{eq:Abound}
\end{align}
Then $|A_0|\le A$. For a damping gain $\kappa_d\ge0$, define
\begin{equation}
 D=\kappa_d\frac{v_+^4}{h_-}c_z,\qquad
 L_T=s_-^2-h_+\left(2g+\frac{\beta w_+^2}{H+h_-}+2A\right),\label{eq:Tbounds}
\end{equation}
and
\begin{align}
 L_u={}&s_-^2-(1-\alpha_0)v_+^2-2gh_+
 -\frac{\beta h_+w_+^2}{H+h_-}-2h_+A-D,\label{eq:ulow}\\
 U_u={}&\left(1+\frac{h_+}{H+h_-}\right)s_+^2
 -(1-\alpha_0)v_-^2-2gh_-+2h_+A+D.\label{eq:uhigh}
\end{align}
\begin{proposition}[Sufficient force feasibility]\label{prop:bounds}
At every state in the tracking region, for the input
\begin{equation}
 \ur=\un+\kappa_d\frac{v_1^4}{h}z,\label{eq:robust}
\end{equation}
one has $L_u\le\ur\le U_u$ and, for every $|d|\le\bar d$,
\begin{equation}
 T_P/\lambda\ge L_T-D-\bar d\,v_+^2.\label{eq:Trob}
\end{equation}
Consequently the inequalities
\begin{equation}
 L_u\ge u_{\mathrm{act},-},\quad U_u\le u_{\mathrm{act},+},\quad
 L_T-D-\bar d\,v_+^2\ge0\label{eq:forcefeas}
\end{equation}
are sufficient for admissible actuation and nonnegative pickup tension throughout that region. For the nominal controller take $\kappa_d=\bar d=0$.
\end{proposition}
\begin{proof}
Using \eqref{eq:acceleration}, rewrite the nominal input as
\begin{equation}
 \un=\left(1+\frac{h}{H+h}\right)S-(1-\alpha_0)v_1^2
 -2gh-\frac{\beta h w^2}{H+h}-2hA_0.\label{eq:ualternative}
\end{equation}
Bounding each term by the enclosing intervals proves \eqref{eq:ulow}--\eqref{eq:uhigh}. Omitting a positive term from the lower bound and a negative term from the upper bound makes these estimates more conservative. The added damping term has magnitude at most $D$. The actual pickup tension equals the nominal-coefficient tension plus this damping term minus $d v_1^2$. Dropping the positive term $hS/(H+h)$ in the nominal tension gives $L_T$, and the remaining two terms are bounded below by $-D-\bar d v_+^2$.
\end{proof}
These sufficient conditions allow feasibility to be checked before a trajectory is applied and ensure that no saturation is needed inside the invariant region. Their conservatism leaves open the possibility of a larger viable region.

\subsection{A constructive reference-rate test}
The bounds constrain how a setpoint change may be requested. For a transition from $r_0$ to $r_1$ over duration $T$, the quintic schedule
\begin{equation}
 r(t)=r_0+(r_1-r_0)(10\sigma^3-15\sigma^4+6\sigma^5),\quad \sigma=t/T,\quad 0\le t\le T,
\end{equation}
extended by constants before and after the transition, is $C^2$. It satisfies
\begin{equation}
 R_1=\frac{15|r_1-r_0|}{8T},\qquad
 R_2=\frac{10|r_1-r_0|}{\sqrt3\,T^2}.\label{eq:quintic}
\end{equation}
Thus a reference governor can test candidate durations using \eqref{eq:positive}--\eqref{eq:upperface} and \eqref{eq:forcefeas}, together with the initial error bound. A duration is accepted only when all inequalities pass. Increasing $T$ reduces the derivative demands. A requested transition may nevertheless remain infeasible if its steady height is inaccessible or its initial state fails the operating-region conditions.

\subsection{How to apply the sufficient inequalities}\label{sec:certificate-use}
The design variables are the gains, an energy cap $\Vbar$, the interval $[w_-,w_+]$, a crown margin, and the reference bounds. The cap bounds the two errors, which combine with the reference bounds to give a range of actual heights and apex speeds. The speed interval encloses the internal state. Once these ranges are known, the force inequalities determine whether the computed command is physically available.

The bound on $q$ is sharper than independently maximizing the two error terms. Write
\[
q-\dot r=-\frac{k_1}{\sqrt{k_0}}(\sqrt{k_0}e)+z/2.
\]
The vector $(\sqrt{k_0}e,z/2)$ has norm at most $\sqrt{2\Vbar}$, so Cauchy--Schwarz gives exactly $c_q$ in~\eqref{eq:cdefs}. Geometrically, this is the support function of the error ellipse in the direction relevant to apex motion.

\begin{figure}[htbp]\centering
\includegraphics[width=\linewidth]{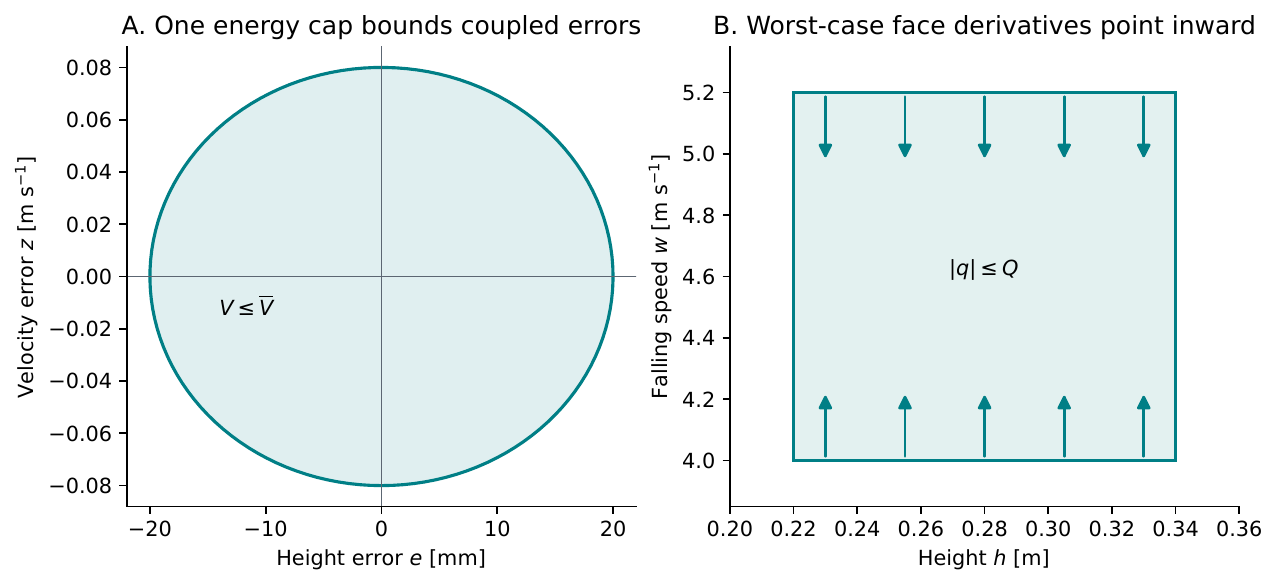}
\caption{Two linked parts of the invariance argument, using the numerical example. The left panel shows the weighted error ellipse. The right panel projects the invariant region onto $(h,w)$, with arrows indicating the inward signs at the falling-speed faces. These arrows do not form a full phase portrait. The actual region also constrains $e,z,q$ and depends on the reference.}\label{fig:tubegeometry}
\end{figure}

At the lower face, the smallest height and largest positive apex speed give the worst case, because increasing $q$ increases $(w+q)^2$ and reduces the falling acceleration. At the upper face, the worst outward derivative occurs at the largest height and $q=-Q$. The condition $w_->2Q$ keeps the relevant crown speeds positive and ensures the monotonicity needed for both evaluations. We thus obtain the two explicit face inequalities and, once these hold, the lower bounds on $v_1=w+2q$ and $s=w+q$ follow.

\begin{figure}[htbp]\centering
\includegraphics[width=.96\linewidth]{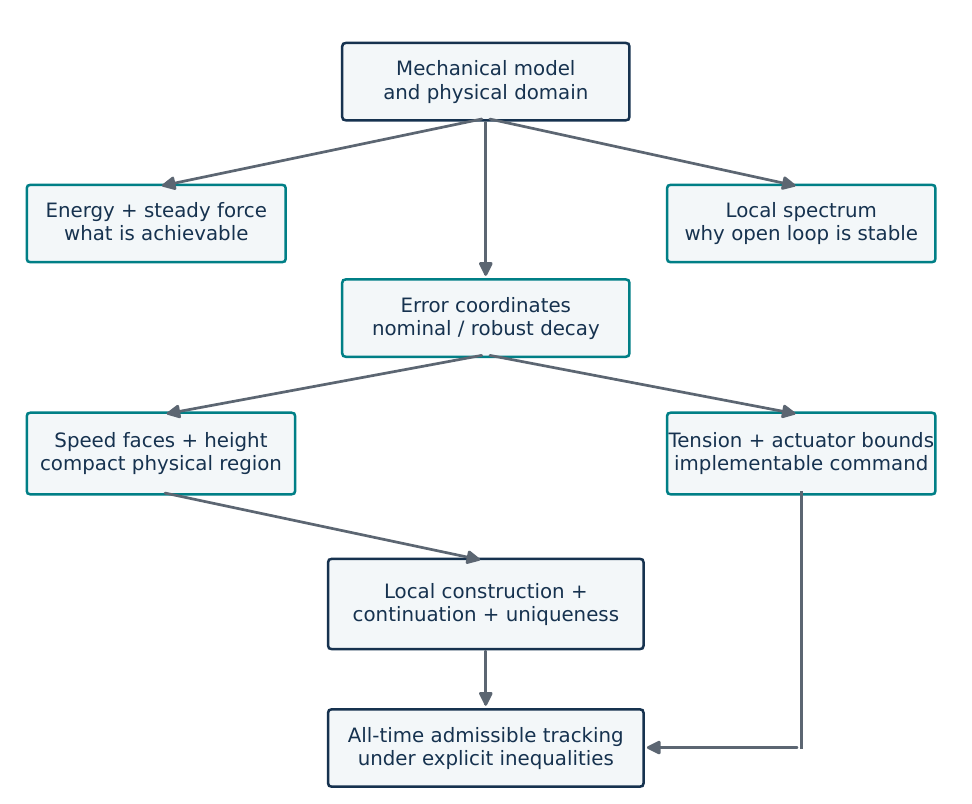}
\caption{Logical route through the deterministic results. Arrows indicate mathematical use in the exposition, rather than the full kernel dependency graph. Force feasibility and the existence result together establish admissible tracking. The local-spectrum result explains the natural operating point but is not needed to derive the nominal error cancellation.}\label{fig:proofroute}
\end{figure}

To establish global existence, we first work on a finite interval where a local solution is already known to exist. The error estimate bounds $h$ and $q$, and the speed comparison then confines the state to a compact box. Local regularity and this compact containment exclude a finite maximal endpoint, so the argument does not presuppose an all-time trajectory. Figure~\ref{fig:proofroute} distinguishes the estimates on finite intervals from the continuation conclusion they support.

\paragraph{A worked calculation.}
For the parameters of Section~\ref{sec:numerics}, $r_-=0.24$, $r_+=0.32$, $R_1=0.02\pi$ and $R_2=0.01\pi^2$. The cap $\Vbar=0.0008$ gives
\[
c_e=0.0200,\quad c_z=0.0800,\quad c_q\simeq0.05657,
\quad Q\simeq0.11940.
\]
Thus $h\in[0.22,0.34]$, $v_1\ge4-2Q\simeq3.76120$, and $s\ge4-Q\simeq3.88060$. At the lower face, the supply $g(H+h_-)=21.7782$ exceeds the required bound of approximately $16.64946$. At the upper face, the supply $22.9554$ lies below the opposing bound of approximately $25.27169$. All these quantities have squared-speed units, and both inequalities have positive slack.

The initial error energy is $0.0001625$, and the robust example has an ultimate cap of $0.000625$, below $\Vbar$. The separate conservative force bounds are approximately $[-19.71183,17.85059]$, within the assumed interval $[-25,25]$, while the normalized pickup-tension lower bound is $6.02325>0$. Evaluated before inspecting the simulated trajectories, these numbers establish a sufficient region. They do not identify the largest possible region.

\paragraph{A reproducible design sequence.}
Begin with a force range and physically meaningful height and speed margins. After choosing a smooth reference family and gains, propose an initial-error cap and speed interval, then evaluate the geometric conditions, the robust cap when needed, and all three force inequalities. The motion is accepted only when every test passes. If a test fails, the appropriate adjustment depends on its cause. A slower reference can reduce derivative demand, but cannot make an inaccessible steady height attainable. Similarly, increasing damping can improve an error bound while making the force condition harder to satisfy. Because these tests are sufficient conditions, the search need not establish feasibility for every physically possible motion.

\FloatBarrier\section{Bounded pickup uncertainty}\label{sec:robust}
Let $G(x)=v_1^2/h$. The gain in \eqref{eq:robust} has units $\mathrm{s^3\,m^{-2}}$, so the added input has the required units. Larger damping reduces sensitivity in the error equation at the possible cost of greater force demand, as the bound $D$ makes explicit.

\begin{theorem}[Disturbance estimate and robust invariant region]\label{thm:robust}
Let $\kappa_d>0$ and Assumption~\ref{ass:parameters} hold. With the unfiltered input \eqref{eq:robust}, on every interval of existence in $\D$,
\begin{align}
 \dot z&=-2k_0e-k_2z-\kappa_dG^2z+Gd,\label{eq:zdrob}\\
 \dot V&\le-2\kappa V-\frac{\kappa_d}{8}G^2z^2
 +\frac{d^2}{8\kappa_d}.\label{eq:Vrob}
\end{align}
In particular,
\begin{equation}
 V(t)\le e^{-2\kappa t}V(0)
 +\frac{\bar d^2}{16\kappa\kappa_d}(1-e^{-2\kappa t}).\label{eq:Vrobsolution}
\end{equation}
If the geometric conditions of Theorem~\ref{thm:tube} hold and
\begin{equation}
 \Vbar\ge\max\left\{V(0),\frac{\bar d^2}{16\kappa\kappa_d}\right\},\label{eq:robcap}
\end{equation}
then \eqref{eq:tube} holds for all time. If \eqref{eq:forcefeas} also holds, the actuator and pickup constraints hold without clipping, for every such disturbance realization.
\end{theorem}
\begin{proof}
The coefficient perturbation adds $Gd$ to $\dot v_1$. Since $\dot a$ depends on $f_2$ and the kinematics, it receives no direct disturbance term. Substitution gives \eqref{eq:zdrob}. Therefore
\[
 \dot V=-k_0k_1e^2-\tfrac14k_2z^2-\tfrac14\kappa_dG^2z^2+\tfrac14Gzd.
\]
Young's inequality gives $Gzd/4\le\kappa_dG^2z^2/8+d^2/(8\kappa_d)$, proving \eqref{eq:Vrob}. Scalar comparison gives \eqref{eq:Vrobsolution}. Under \eqref{eq:robcap}, the level $V=\Vbar$ cannot be crossed outward. The speed-face argument of Theorem~\ref{thm:tube} uses only the bound on $V$, the kinematics, and the falling-leg equation, so it applies unchanged. Proposition~\ref{prop:global} justifies the local Carath\'eodory construction and global continuation using these estimates on finite existence intervals. Proposition~\ref{prop:bounds} then establishes the force constraints.
\end{proof}
For the dimensionally weighted error $\xi=(\sqrt{k_0}e,z/2)$, $\norm{\xi}^2=2V$. Integration yields
\begin{equation}
 \int_0^T\norm{\xi(t)}^2\dd t
 \le\frac{V(0)}{\kappa}+\frac{1}{8\kappa\kappa_d}\int_0^T d(t)^2\dd t.\label{eq:L2}
\end{equation}
The corresponding zero-initial-error $L_2$ gain is at most $(8\kappa\kappa_d)^{-1/2}$, on trajectories whose existence has been established. For persistent bounded uncertainty, the $L_2$ energy may be infinite. In that case \eqref{eq:Vrobsolution} gives the appropriate ultimate bound.

\subsection{What robustness costs}
The disturbance is multiplicative in the physical model and enters the error equation through the known factor $G=v_1^2/h$. We choose the damping term $-\kappa_dG^2z$ so that the cross term generated by $Gd$ in the energy derivative can be bounded without first replacing $G$ by a global constant. The required elementary inequality is
\[
0\le\left(\sqrt{\kappa_d}Gz-\frac{d}{\sqrt{\kappa_d}}\right)^2,
\qquad
\frac14Gzd\le\frac18\kappa_dG^2z^2+\frac{d^2}{8\kappa_d}.
\]
One half of the added damping absorbs the cross term, leaving the other half as additional dissipation in~\eqref{eq:Vrob}. Omitting this favorable term gives the scalar bound. Including it can sharpen particular estimates, although it is unnecessary for the uniform error bound.

The asymptotic height bound and the region-wide additional force allowance scale differently:
\[
\limsup_{t\to\infty}|h(t)-r(t)|\le
\frac{\bar d}{\sqrt{8\kappa\kappa_d k_0}},
\qquad
D=\kappa_d\frac{v_+^4}{h_-}c_z.
\]
The error bound decreases like $\kappa_d^{-1/2}$, while the force allowance increases linearly. Any increase in the gain must therefore be checked against actuator capacity and pickup tension. The functions $r,d$ need not be periodic, and the theorem requires no probability law for $d$.

\begin{figure}[htbp]\centering
\includegraphics[width=\linewidth]{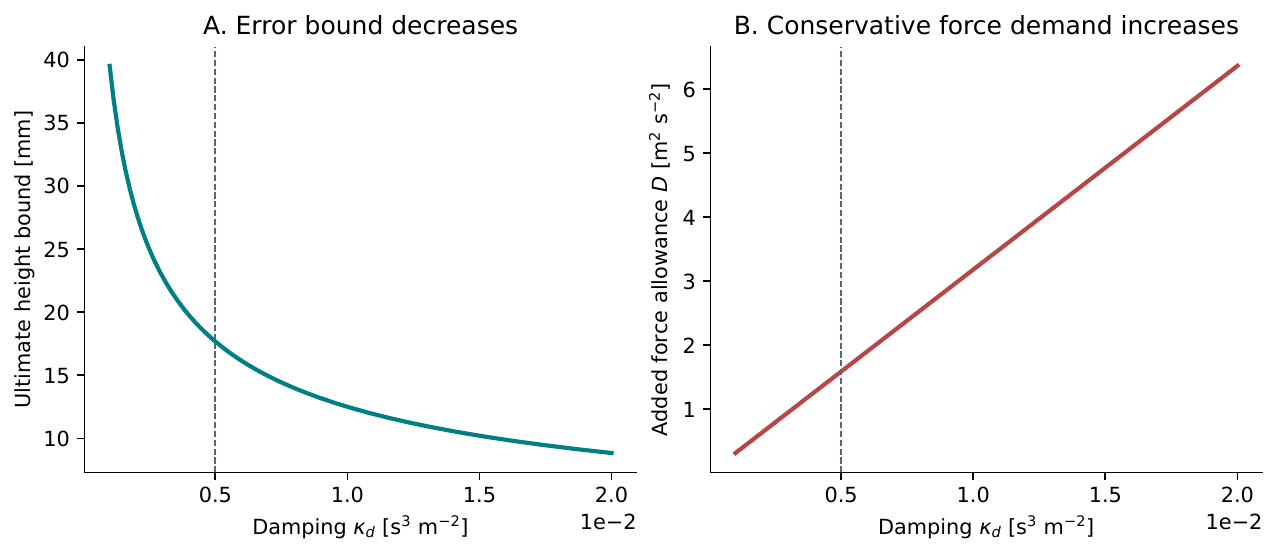}
\caption{Analytical damping tradeoff with the reference, gains and proposed region of the numerical example held fixed. The dashed line marks $\kappa_d=0.005$. Parts of the displayed gain range may fail the robust cap or force tests, so the curves alone do not establish admissibility for every gain.}\label{fig:tradeoff}
\end{figure}

The finite-energy estimate answers a different question from the ultimate bound. For $d\in L^2(0,\infty)$ and zero initial error, it bounds the ratio of accumulated weighted error energy to accumulated disturbance energy. A persistent disturbance can have infinite $L^2$ norm while remaining uniformly bounded, in which case the pointwise envelope is still meaningful. The deterministic theorem also covers a random bounded coefficient path by path whenever the stated bound holds, without requiring independence, stationarity or a Gaussian hypothesis.

\FloatBarrier\subsection{Global construction and continuation}\label{sec:global}
We now use the preceding estimates to construct the solution from its initial data and continue it for all forward time. This connects the tracking bounds to an existence theorem within the same reduced model, without adding a mechanical assumption.

\begin{proposition}[Global existence, uniqueness, and physical tracking]\label{prop:global}
Suppose the reference satisfies \eqref{eq:refbounds}, the geometric inequalities \eqref{eq:positive}--\eqref{eq:upperface} hold, and $w(0)\in[w_-,w_+]$. Consider either the nominal law, with $d\equiv0$, $\kappa_d=0$, and $V(0)\le\Vbar$, or the robust law, with $\kappa_d>0$, measurable $|d(t)|\le\bar d$, and \eqref{eq:robcap}. Then the initial-value problem has a unique solution on $[0,\infty)$, locally absolutely continuous and satisfying the controlled equations almost everywhere. It remains in
\begin{equation}
 K=[h_-,h_+]\times[v_-,v_+]\times[w_-,w_+],\qquad
 v_-=w_--2Q>0,\quad v_+=w_++2Q.\label{eq:compactK}
\end{equation}
The nominal or robust tracking estimate, as appropriate, and the crown margin hold for every $t\ge0$. If \eqref{eq:forcefeas} holds, the unfiltered command also respects the actuator and pickup constraints. For continuous $d$, including $d\equiv0$, the solution is classical. These statements concern the ideal sustained-supply model specified in Section~\ref{sec:model}.
\end{proposition}
\begin{proof}
Let $U=\{(h,v_1,w):h>0\}$. The closed-loop field can be written as $F(t,x)=\Phi(p(t),x)$, where $p=(r,\dot r,\ddot r,d)$ and $\Phi$ is smooth in its arguments on $h>0$. Indeed, its only state denominators are powers of $h$ and $H+h$. The reference bounds and $|d|\le\bar d$ place $p$ in a fixed compact box. On each compact state box contained in $U$, both $\Phi$ and its state derivative are therefore bounded uniformly in $p$. Thus $F$ is measurable in time, continuous and locally Lipschitz in the state, with locally integrable bounds. The local existence and uniqueness theorem for Carath\'eodory equations therefore applies. For continuous $d$, it reduces to the classical Picard--Lindel\"of theorem \cite[Sections 2.2--2.3]{teschl}. Denote the forward maximal interval by $[0,\tau)$.

Fix $T<\tau$. The error-coordinate chain rule is valid almost everywhere because the trajectory is absolutely continuous and the coordinate transformation is smooth in the state. The nominal Lyapunov identity, or \eqref{eq:Vrob} in the robust case, can therefore be integrated after multiplication by $e^{2\kappa t}$. It yields the corresponding estimate for $V$ on $[0,T]$ and hence $V\le\Vbar$. This step does not assume a bound on either transport velocity. The coordinate estimates then give $h_-\le h\le h_+$ and $|q|\le Q$.

For completeness, the weak inequalities at the speed faces suffice even when a trajectory touches a face. Write the right side of \eqref{eq:internal} as $f(t,w)$, using the already bounded functions $h(t)$ and $q(t)$. For either constant face $b\in\{w_-,w_+\}$,
\begin{equation}
 f(t,w)-f(t,b)=a_b(t)(w-b),\qquad
 a_b(t)=-\frac{(1-\beta)(w+b)+2q}{H+h}.\label{eq:facefactor}
\end{equation}
The face inequalities give $f(t,w_-)\ge0$ and $f(t,w_+)\le0$. Consequently $y_-=w_--w$ and $y_+=w-w_+$ satisfy $\dot y_\pm\le a_{w_\pm}(t)y_\pm$ almost everywhere. On $[0,T]$ the coefficients are continuous, and the integrating factors preserve $y_\pm\le0$. It follows that $w\in[w_-,w_+]$, $v_1=w+2q\in[v_-,v_+]$, and $s=w+q\ge w_--Q\ge v_{\min}$. Since $T<\tau$ was arbitrary, $x(t)\in K$ for every $t<\tau$.

Suppose now that $\tau<\infty$. The compact set $K$ lies strictly inside $U$. A useful explicit bound on the vector field is obtained by setting
\begin{align}
 G_+&=\frac{v_+^2}{h_-},&
 M_w&=g+\frac{s_+^2+\beta w_+^2}{H+h_-},\nonumber\\
 M_v&=M_w+2A+\kappa_dG_+^2c_z+G_+\bar d,&
 M&=\max\{Q,M_w,M_v\}.\label{eq:continuationM}
\end{align}
Here $A,c_z,s_+$ are the previously defined bounds, and the nominal case uses $\kappa_d=\bar d=0$. The falling-leg equation gives $|\dot w|\le M_w$. The identity
\[
 \dot q=A_0-\tfrac12\kappa_dG^2z+\tfrac12Gd
\]
and $\dot v_1=\dot w+2\dot q$ give $|\dot v_1|\le M_v$, while $|\dot h|\le Q$. Thus $\|\dot x\|_\infty\le M$ almost everywhere and
\begin{equation}
 \|x(t)-x(s)\|_\infty\le M|t-s|,\qquad 0\le s,t<\tau.\label{eq:timeLipschitz}
\end{equation}
The trajectory has a limit $x_*\in K$ as $t\uparrow\tau$. Apply the same local existence theorem at $(\tau,x_*)$ and concatenate the integral solutions. Because the solutions share their endpoint, their concatenation is absolutely continuous and satisfies the integral equation. This construction requires no pointwise derivative at a jump of $d$. This extends the solution beyond $\tau$, a contradiction. Therefore $\tau=\infty$. The uniform state Lipschitz bound on compact sets also gives uniqueness by Gronwall's inequality. The force and tension statements follow from Proposition~\ref{prop:bounds} along the solution just constructed.
\end{proof}

\paragraph{Formal realization.}
The companion proves both classical and measurable-input constructions, the general AC/integral equivalence, and uniqueness among forward locally absolutely continuous competitors. Appendix~\ref{app:formal} identifies the entry points. Both the written continuation argument and the formal bounded-extension construction obtain the global trajectory from estimates on finite existence intervals.

\section{What a crown barrier can and cannot guarantee}\label{sec:barrier}
The operating-region result gives joint physical guarantees under explicit conditions. We next examine a crown barrier as an additional local mechanism and determine exactly which constraints it protects. Its scope is narrower than invariance of the full domain $\D$. Let
\begin{equation}
 B_C(x)=s-v_{\min},\qquad \gamma>0.
\end{equation}
For bounded $\alpha\in[\amin,\amax]$, the conditions $T_P\ge0$ and $\dot B_C+\gamma B_C\ge0$ hold for all coefficient values if and only if
\begin{align}
 u&\ge L_B:=-(1-\amax)v_1^2,\label{eq:barlower}\\
 u&\le U_B:=S-(1-\amin)v_1^2-gh+h[f_2(x)+2\gamma B_C(x)].\label{eq:barupper}
\end{align}
The worst cases differ: large $\alpha$ lowers the pickup tension, while small $\alpha$ lowers the rising acceleration for fixed input.

\begin{proposition}[Exact pointwise barrier feasibility]\label{prop:barrier}
Before actuator limits are imposed, the robust barrier interval is nonempty exactly when
\begin{equation}
 \frac{HS+\beta h w^2}{H+h}+2h\gamma B_C
 \ge(\amax-\amin)v_1^2.\label{eq:barwidth}
\end{equation}
With actuator limits, set
\begin{equation}
 L=\max\{u_{\mathrm{act},-},L_B\},\qquad
 U=\min\{u_{\mathrm{act},+},U_B\}.\label{eq:combinedinterval}
\end{equation}
The least-change filter exists precisely when $L\le U$, and then is
\begin{equation}
 u^*=\min\{U,\max\{L,u_{\mathrm{des}}\}\}.\label{eq:projection}
\end{equation}
Along any solution remaining in the model domain on which this feasible input is applied, $B_C(0)\ge0$ implies $B_C(t)\ge0$, and $T_P(t)\ge0$.
\end{proposition}
\begin{proof}
Differentiating $B_C$ gives $(f_1+f_2)/2$. Solving its affine input inequality and the pickup inequality gives \eqref{eq:barlower}--\eqref{eq:barupper}. Their difference simplifies to the left side of \eqref{eq:barwidth} minus its right side. Intersecting with the actuator interval proves feasibility and the scalar projection formula. Finally $\dot B_C\ge-\gamma B_C$ gives $B_C(t)\ge B_C(0)e^{-\gamma t}$ by comparison. This last argument does not establish that the solution stays in $\D$.
\end{proof}
In the nominal case $\amax=\amin$, the unconstrained interval is nonempty whenever $B_C\ge0$. Under uncertainty, the interval can be empty even when the height and both velocities are positive, as the example in Appendix~\ref{app:counter} shows. The constraints then conflict, and a clipping function called with reversed limits cannot supply a feasible input.

This is a CBF-based least-change filter in the terminology of \cite{ames}. It is not, by itself, a combined CLF--CBF guarantee. If the desired input is \eqref{eq:robust}, applying a different $u$ adds
\begin{equation}
 -\frac{z}{4h}(u-\ur)\label{eq:clfchange}
\end{equation}
to $\dot V$. To preserve \eqref{eq:Vrob}, one can require this term to be nonpositive. The resulting affine input constraint must be checked together with the barriers, since the conditions may conflict. A relaxation can keep a numerical program solvable, but changes the tracking estimate. The theorems of Sections~\ref{sec:tube}--\ref{sec:robust} use the unfiltered law with proved force bounds, and do not assume that projection preserves cancellation.

\subsection{An interval calculation, and its exact logical scope}
Because the input is scalar, each affine inequality defines a half-line of admissible commands. Pickup tension imposes a lower bound on the input, while the crown derivative imposes an upper bound. Their intersection is nonempty exactly when the lower endpoint does not exceed the upper endpoint. Subtracting~\eqref{eq:barlower} from~\eqref{eq:barupper} gives
\begin{align*}
U_B-L_B
&=S+h(f_2-g)+2h\gamma B_C-(\alpha_+-\alpha_-)v_1^2\\
&=\frac{HS+\beta hw^2}{H+h}+2h\gamma B_C
 -(\alpha_+-\alpha_-)v_1^2.
\end{align*}
The compatibility condition is therefore that this difference be nonnegative. Uncertainty reduces the available width by $(\alpha_+-\alpha_-)v_1^2$. After intersecting with the actuator interval, the input closest to the requested command is its Euclidean projection onto the resulting interval.

\begin{figure}[htbp]\centering
\includegraphics[width=\linewidth]{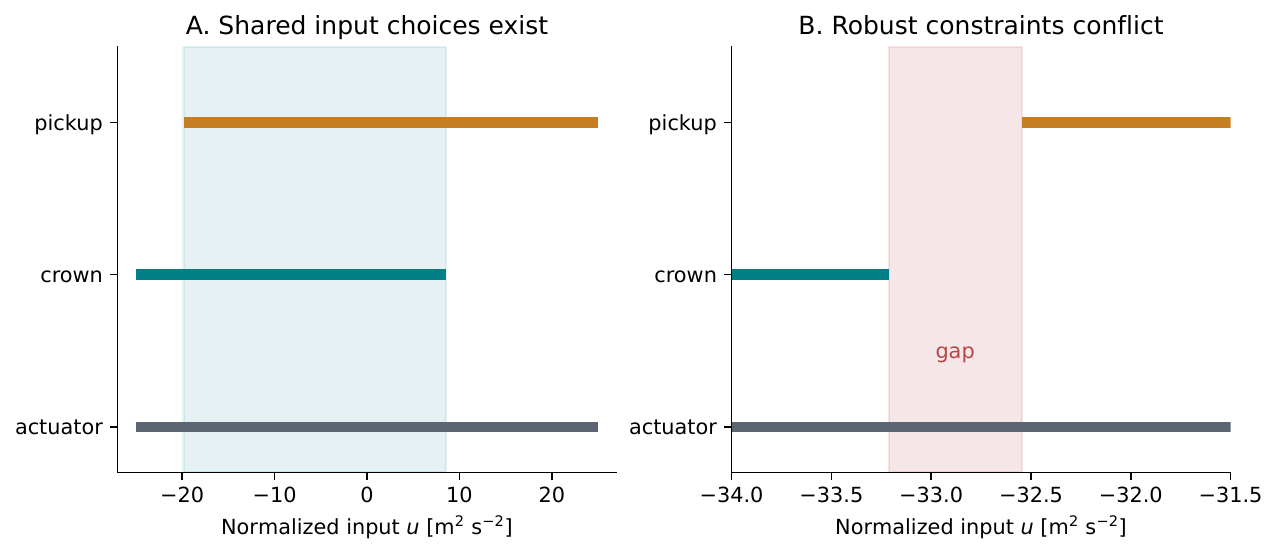}
\caption{Input-space geometry. Left: nominal equilibrium parameters with coefficient interval $[0.11,0.13]$, $\gamma=10$ and actuator bounds $[-25,25]$. Right: the counterexample in Appendix~\ref{app:counter}, shown with a wider illustrative actuator interval $[-34,-31.5]$ to expose the incompatibility of the two barrier constraints themselves. Shading in the right panel is a gap, not an admissible set.}\label{fig:intervals}
\end{figure}

Multiplying $\dot B_C\ge-\gamma B_C$ by the integrating factor $e^{\gamma t}$ yields $B_C(t)\ge e^{-\gamma t}B_C(0)$. This proves the crown claim on intervals where the model and feasible input apply, but does not by itself address $h=0$ or $w=0$. Projection also changes the tracking equation. The tracking estimate is preserved when the projection is inactive or satisfies the additional sign condition in~\eqref{eq:clfchange}. The unfiltered region theorem and the conditional barrier theorem must therefore be applied with their respective hypotheses.

\FloatBarrier\section{White noise is a separate model}\label{sec:stochastic}
\subsection{From bounded fluctuations to a Brownian forcing hypothesis}
The bounded-input theorem treats $d(t)$ as an ordinary measurable function. A diffusion approximation instead postulates
\begin{equation}
 \dd v_1=f_1(x,u,\alpha_0)\dd t+\sigma G(x)\dd W_t,
 \qquad G(x)=v_1^2/h,\label{eq:diffusion}
\end{equation}
with no diffusion in $h$ or $w$. Here $W$ is Brownian motion. For correspondence with the formal construction we use its natural filtration for the integrals and the augmented filtration for stopping times, and $\sigma$ has units $\mathrm{s}^{1/2}$. The expression $d(t)=\sigma\dot W_t$ is a mnemonic for an integral equation involving the It\^o integral, since Brownian motion has no ordinary time derivative. Consequently, a white-noise coefficient has no pointwise bound comparable to $|d|\le\bar d$.

Introducing this approximation makes a separate modeling assumption beyond the earlier uncertainty bound. Relating it to a physical process with finite correlation would require an appropriate limiting argument and calibration. Here the model lets us examine which guarantees change when a diffusion component is introduced.

\subsection{Constructing a solution before applying an energy formula}
Under the unfiltered damped controller and a $C^2$ reference with bounded $r$, $\dot r$ and $\ddot r$ as in~\eqref{eq:refbounds}, the mechanical drift and diffusion are continuous in time and locally Lipschitz in state on $h>0$. A nested sequence of compact positive-height regions admits globally Lipschitz coefficient extensions. Each extended system has a continuous adapted strong solution, and compatible solutions agree until their common exit. Patching these solutions gives a maximal local process with a positive stopping lifetime $\zeta$. A finite lifetime can occur only if the process leaves every compact subset of the regular domain. The construction establishes existence and compatibility up to that lifetime, without proving $\zeta=\infty$.

The Lean development follows this route using actual predictable stochastic integrals. It identifies their continuous versions simultaneously in time, proves genuine exit-time properties, and transfers the localized equations to the maximal process. Uniqueness is established within the constructed strong-solution/cutoff interface, which suffices for the results below without an unrestricted weak-solution uniqueness claim. For a physical interpretation, the process must also be stopped before any relevant transport boundary or failure of the model assumptions.

\subsection{The energy correction is unavoidable}
The variables $h,w$ have finite-variation trajectories before localization, although their derivatives are random through the state. Consequently $a=w+2\dot r-2k_1e$ has no Brownian component. The error equations are
\begin{align*}
\dd e&=(-k_1e+z/2)\dd t,\\
\dd z&=(-2k_0e-k_2z-\kappa_dG^2z)\dd t+\sigma G\dd W_t.
\end{align*}
Apply It\^o's formula to the Lyapunov function $V=k_0e^2/2+z^2/8$. The second derivative in $z$ is $1/4$, so the half-Hessian term is $\sigma^2G^2/8$. Thus
\begin{equation}
 \mathcal LV=-k_0k_1e^2-\tfrac14k_2z^2
 -\tfrac14\kappa_dG^2z^2+\tfrac18\sigma^2G^2.\label{eq:generator}
\end{equation}
More explicitly,
\[
\dd V=\mathcal LV\dd t+\frac{\sigma zG}{4}\dd W_t.
\]
The positive correction comes from quadratic variation. It would be absent if one tried to apply the deterministic Young inequality by assigning white noise a fictitious pointwise value.

\begin{proposition}[Localized stochastic energy estimate]\label{prop:stochastic-energy}
Take a compact localization region on which $|G|\le G_{\max}$, and let $\tau$ be its genuine first exit for the constructed mechanical solution. On each finite deterministic horizon on which the reference bounds hold, put $\theta=t\wedge\tau$. Then
\begin{equation}
 \mathbb E V(\theta)+2\kappa\mathbb E\int_0^\theta V(s)\dd s
 \le V(0)+\frac{\sigma^2G_{\max}^2}{8}\mathbb E\theta.\label{eq:stopped}
\end{equation}
For the compact exhaustion, the same statement is applied with $\theta=t\wedge\tau_n\wedge T$, where $T<\infty$. If the $n$th compact level gives $\|x\|\le\rho_n$ and $h\ge1/\rho_n$, one may take $G_{\max}=\rho_n^3$.
\end{proposition}
\begin{proof}
Drop the additional negative damping term from~\eqref{eq:generator} and use $k_1,k_2\ge\kappa$. The resulting bound is $\mathcal LV\le-2\kappa V+\sigma^2G_{\max}^2/8$. On the stopped finite horizon, the martingale integrand $\sigma zG/4$ is bounded by compact localization and the reference bounds. Its squared integral has finite expectation. The stopped It\^o integral is therefore a square-integrable martingale of mean zero. Integrate the energy identity to $\theta$, take expectations, and apply the generator bound. The constructed process satisfies the original coefficients until the exit, so the estimate transfers from the cutoff realization to the maximal solution.
\end{proof}

Both the terminal energy and the integral on the left are evaluated up to the same random time. An integral extending to deterministic $t$ would require a different statement. The estimate therefore describes energy before exit, leaving global stochastic existence, an invariant probability distribution and almost-sure crown safety to be addressed separately.

\subsection{The actual crown diffusion and the small-time obstruction}
Since $B_C=(v_1+w)/2-v_{\min}$ is linear in the velocities,
\begin{equation}
 \dd B_C=\tfrac12(f_1+f_2)\dd t+\tfrac12\sigma G\dd W_t.\label{eq:barriernoise}
\end{equation}
The coefficient normal to the crown boundary is therefore
\[
a(x)=\frac{\sigma v_1^2}{2h}.
\]
Stochastic viability has long required a joint analysis of drift and diffusion. Aubin and Da Prato develop this through stochastic tangent conditions~\cite{aubin1990}, and Milian gives viability conditions for polyhedra together with a comparison theorem~\cite{milian1995}. The crown constraint is locally an affine half-space. Its normal diffusion is exactly the coefficient above, which is nonzero at a positive-transport boundary state when $\sigma\ne0$. This is the obstruction that an inward drift condition alone cannot remove.

Here we make the obstruction explicit on the constructed mechanical process. The Brownian law of the iterated logarithm and the continuous-martingale time-change theorem are classical background~\cite[Sections 2.3 and 5.8]{baudoin}. Our proof uses local freezing of the actual integral, proved below without assuming a time-change representation. It then identifies the sharp positive and negative small-time excursions before a genuine local exit. The conclusion is thus tied to the coefficient of the mechanical SDE and its solution, rather than to a separately prescribed Brownian trajectory.

\begin{theorem}[Small-time freezing of the stochastic integral]\label{thm:freezing}
Let $a_t$ be a continuous scalar coefficient adapted to the natural filtration of Brownian motion $W$ starting at zero, with deterministic initial value $a_0$. We use fixed continuous versions of $a$ and $W$. No global deterministic bound on $a$ is assumed. Define the actual localized It\^o integral $N_t=\int_0^t a_s\,\dd W_s$. With
\[
\ell(t)=\sqrt{2t\log\log(1/t)},\qquad 0<t<e^{-1},
\]
one has, almost surely,
\[
N_t=a_0W_t+o(\ell(t))\quad(t\downarrow0).
\]
The conclusion is unaffected by restriction to any almost surely positive initial stopping interval. Time in the logarithm is expressed in a fixed unit; changing that unit does not change the limiting statement.
\end{theorem}
\begin{proof}
First clip $a$ at the deterministic level $|a_0|+1$. Continuity at zero and prefix locality of the actual stochastic integral give an almost surely positive initial interval on which the clipped and original integrals agree. It therefore suffices to prove the statement for this bounded coefficient. The continuous-adapted and predictable localized integrals used here are identified simultaneously in time through their common Riemann-sum limits and continuous versions.

Linearity reduces the bounded case to the coefficient $b_t=a_t-a_0$, continuous with $b_0=0$. If a predictable coefficient $c$ is bounded by $\varepsilon$, its stochastic integral $M$ has the exponential estimate
\[
\mathbb P\left(\sup_{0\le s\le T}|M_s|\ge x\right)
\le 2\exp\left(-\frac{x^2}{2\varepsilon^2T}\right),\qquad x,T>0.
\]
One obtains it from the exponential martingale and a maximal inequality, separately for $M$ and $-M$. On geometric times $T_n=\vartheta^n$, $0<\vartheta<1$, choose $x_n=C\varepsilon\sqrt{T_n\log n}$ with $C$ large enough to make the probabilities summable. Borel--Cantelli gives an almost-sure eventual running-maximum bound. Comparing $t\in[T_{n+1},T_n]$ with the two endpoints gives
\[
\limsup_{t\downarrow0}\frac{|M_t|}{\ell(t)}\le C_\vartheta\varepsilon
\]
for a fixed finite constant. Its sharp value is immaterial.

Apply this estimate to $b$ clipped to $[-\varepsilon,\varepsilon]$. Continuity at zero makes clipping inactive on a positive initial interval of each continuous path. Prefix locality of the actual stochastic integral then identifies the clipped and unclipped integrals on that interval. Taking the intersection of the conull events for the countable choices $\varepsilon=1/(n+1)$ and letting $n\to\infty$ shows that the normalized residual tends to zero. The checked proof constructs continuous integral versions and establishes their prefix locality directly, without assuming a random time-change representation.
\end{proof}

\begin{theorem}[Mechanical crown LIL and boundary non-invariance]\label{thm:mechanical-lil}
Consider the constructed maximal local mechanical solution with deterministic initial state $h_0>0$, $v_{1,0}>0$, $w_0>0$, and $\sigma\ne0$. Let $a_0=\sigma v_{1,0}^2/(2h_0)$ and $B_C(0)=0$. Almost surely,
\begin{equation}
\limsup_{t\downarrow0}\frac{B_C(t)}{|a_0|\ell(t)}=1,
\qquad
\liminf_{t\downarrow0}\frac{B_C(t)}{|a_0|\ell(t)}=-1.\label{eq:mechanical-lil}
\end{equation}
For every $\varepsilon>0$ there is a time $0<t<\varepsilon$, before the positive genuine initial cutoff exit and before the maximal lifetime, for which $B_C(t)<0$. Hence the closed crown constraint is not almost surely invariant from such a boundary state under this diffusion model and the stated locally regular feedback.
\end{theorem}
\begin{proof}
Before a positive compact exit, the crown drift is locally bounded along the continuous path, so its integral $A_t$ satisfies $|A_t|\le Kt$. The mechanical coordinate equations give
\[
B_C(t)=A_t+\int_0^t\frac{\sigma v_1(s)^2}{2h(s)}\,\dd W_s
=a_0W_t+o(\ell(t))+O(t).
\]
The coefficient is continuous at zero and $a_0\ne0$, so Theorem~\ref{thm:freezing} applies. The two-sided Brownian law of the iterated logarithm gives normalized upper and lower limits $1$ and $-1$. Since $t/\ell(t)\to0$, the drift does not change them. If $a_0<0$, the upper and lower limits exchange roles before division by $|a_0|$, giving the same pair. A negative lower limit supplies negative values arbitrarily near zero. The exit and lifetime are positive, so these times may be chosen before both. The formal proof first establishes the coordinate expansion on a measure-preserving continuous realization and then returns it to the original process on one conull event.
\end{proof}

\begin{figure}[htbp]\centering
\includegraphics[width=\linewidth]{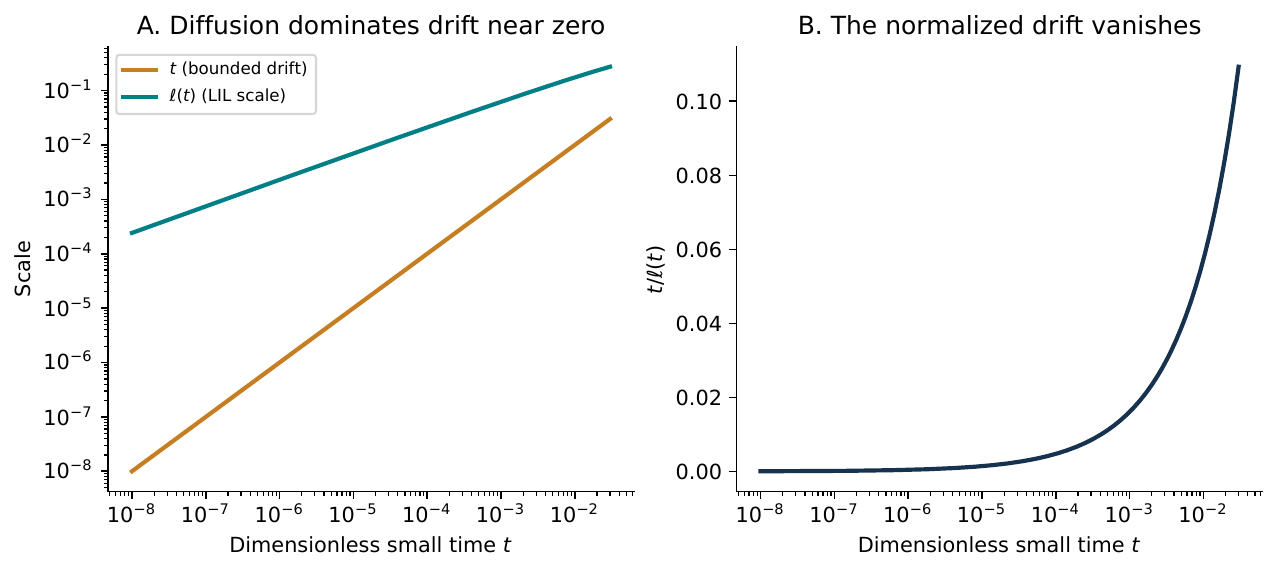}
\caption{The scale comparison behind the boundary theorem. These are deterministic functions, not sampled Brownian paths. Time is nondimensionalized before taking logarithms. The limit is read toward the left.}\label{fig:lilscales}
\end{figure}

\begin{figure}[htbp]\centering
\includegraphics[width=.97\linewidth]{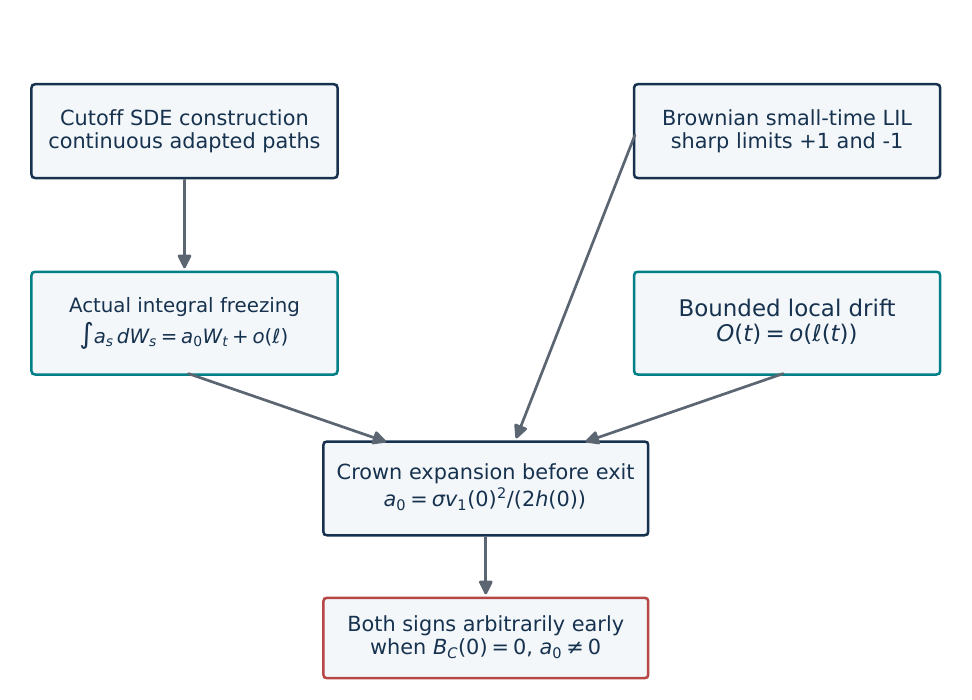}
\caption{The stochastic proof route. Freezing is proved for the actual variable-coefficient integral and transferred to the constructed mechanical process before its genuine exit. This closes the link between a Brownian germ theorem and the mechanical crown variable.}\label{fig:stochasticroute}
\end{figure}

The theorem concerns a process starting on the boundary. It gives no immediate-exit conclusion or exit probability for an interior state. The assumptions exclude $\sigma=0$ as well as singular, reflecting, impulse or diffusion-modifying controls. Within this scope, the result explains why a deterministic drift barrier does not by itself establish stochastic safety. Further questions include finite-horizon chance constraints and models with bounded colored disturbances or a diffusion tangent to the safety boundary. Each of these objectives and models requires its own analysis.

\FloatBarrier\section{Reproducible numerical illustrations}\label{sec:numerics}
We used Python to generate all numerical results in this section. They illustrate \eqref{eq:model} using numerical integration of the reduced equations. Experimental measurements and independent bead-chain simulations are outside this comparison. We choose illustrative parameters
\[
 H=2,\quad g=9.81,\quad\alpha_0=0.12,\quad\beta=0.02,
 \quad k_0=4,\quad k_1=2,\quad k_2=3,
\]
in consistent SI units. These values are not described as a fitted experimental dataset. The reference is
\begin{equation}
 r(t)=0.28+0.04\sin(\pi t/2)\ \mathrm m,\qquad 0\le t\le8\ \mathrm s.\label{eq:numref}
\end{equation}
The initial conditions are specified through $e(0)=0.005\,\mathrm m$, $z(0)=0.03\,\mathrm{m\,s^{-1}}$, $w(0)=\sqrt{g(H+r(0))/(1-\beta)}$, and the inverse of \eqref{eq:errors}. The operating region uses $\Vbar=0.0008\,\mathrm{m^2\,s^{-2}}$, $w_-=4$, $w_+=5.2\,\mathrm{m\,s^{-1}}$, and $v_{\min}=3.2\,\mathrm{m\,s^{-1}}$. For the uncertain case, $\bar d=0.01$ and $\kappa_d=0.005\,\mathrm{s^3\,m^{-2}}$. We specify the normalized actuator interval as $[-25,25]\,\mathrm{m^2\,s^{-2}}$ for this model calculation, without attributing it to a measured device.

\begin{table}[tb]
\centering\small
\begin{tabular}{lrr}
\toprule
Quantity & Value & Unit \\
\midrule
Natural height & 0.279070 & m \\
Natural speed & 4.776395 & m\,s$^{-1}$ \\
Hurwitz determinant & 1164.784792 & s$^{-3}$ \\
Tube: minimum height & 0.220000 & m \\
Tube: maximum height & 0.340000 & m \\
Tube: crown speed lower bound & 3.880600 & m\,s$^{-1}$ \\
Robust pickup tension lower bound & 6.023250 & m$^2$\,s$^{-2}$ \\
Robust control lower bound & -19.711832 & m$^2$\,s$^{-2}$ \\
Robust control upper bound & 17.850592 & m$^2$\,s$^{-2}$ \\
\bottomrule
\end{tabular}

\caption{Computed equilibrium quantities and conservative analytical bounds. Pickup tension is divided by $\lambda$. The control bounds include damping and are inside the assumed actuator interval.}\label{tab:certificate}
\end{table}

We first check the characteristic coefficients against numerical determinants and eigenvalues at 15,251 parameter pairs in the stated triangle, check the Jacobian by central differences, verify the transfer-function numerator independently, and evaluate both sides of the energy identity at 1,000 sampled states. These checks can expose implementation errors, while the analytical proofs establish the statements throughout the stated domains.

For the chosen reference and region, the lower and upper falling-speed face inequalities, positive-transport conditions, force bounds, and robust energy cap all pass. In particular $\bar d^2/(16\kappa\kappa_d)=0.000625<\Vbar$. Table~\ref{tab:certificate} reports bounds obtained directly from the formulas rather than fitted to the trajectories. The ultimate height-error bound from \eqref{eq:Vrobsolution} is approximately $17.68\,\mathrm{mm}$. This conservative analytical bound measures a different quantity from the observed RMS error.

\begin{figure}[tb]
\centering\includegraphics[width=\linewidth]{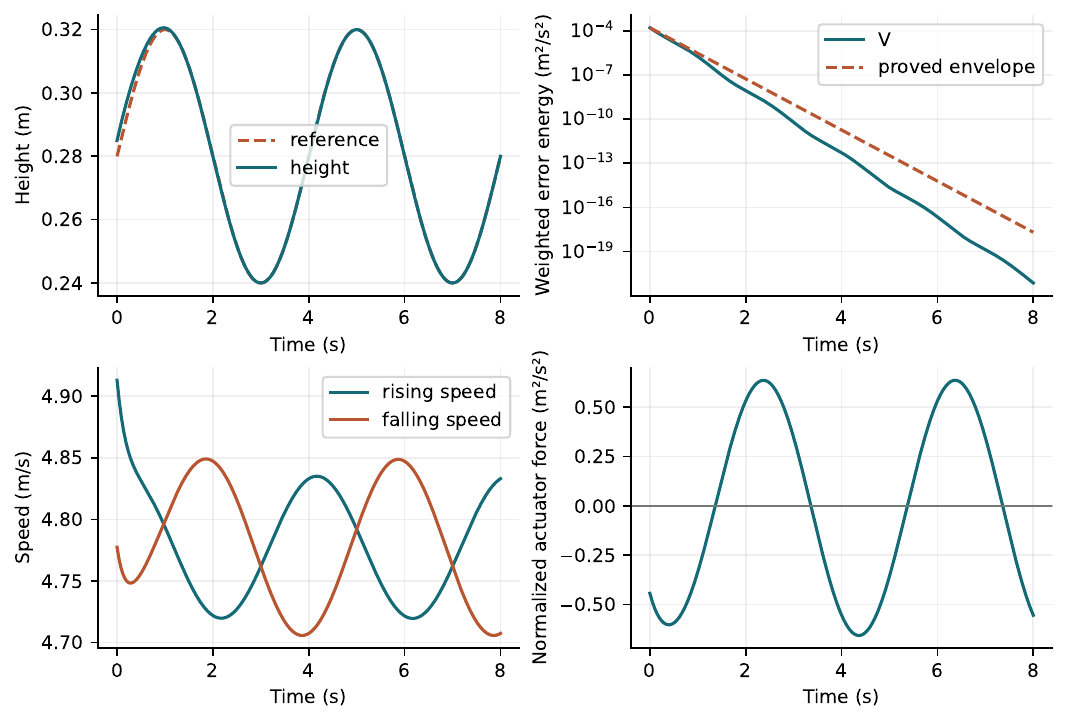}
\caption{Nominal model cancellation for the smooth reference. The weighted error energy stays below its proved exponential envelope. The input changes sign as the desired height passes through the natural height.}\label{fig:tracking}
\end{figure}

Nominal integration uses the explicit DOP853 method with relative tolerance $2\times10^{-9}$, absolute tolerance $2\times10^{-11}$, and maximum step $0.002\,\mathrm s$. Repeating this trajectory with tenfold tighter tolerances and half the maximum step gives a maximum componentwise state difference of \texttt{\RefinementError} in the corresponding SI units. We also compare the simulated errors directly with the matrix-exponential solution of \eqref{eq:errdyn} and test the numerical energy against its envelope. Figure~\ref{fig:tracking} shows the resulting model tracking, whose agreement with an experiment remains to be assessed.

\subsection{Paired bounded-random-input trials}
Twenty paired realizations use seeds $0,\ldots,19$ of NumPy's default random generator. On a grid of spacing $\Delta=0.002\,\mathrm s$, a stationary unit-variance Gaussian sequence is generated by
\begin{equation}
 Y_{j+1}=aY_j+\sqrt{1-a^2}\,\eta_j,\quad a=e^{-\Delta/0.08},\quad
 Y_0,\eta_j\ \hbox{independent standard normal},
\end{equation}
and the coefficient values are $d_j=0.01\tanh Y_j$, interpolated linearly between grid points. Thus $|d(t)|<0.01$ for every realization. The plant is an ODE with a continuous bounded random input, not \eqref{eq:diffusion}. Each seed is applied to both the nominal controller and its damped version. RMS height error and RMS normalized input are measured over $2\le t\le8\,\mathrm s$, on the same output grid.

\begin{table}[tb]
\centering\small\begin{tabular}{lrr}
\toprule
Metric (20 paired realizations) & Nominal & Damped \\
\midrule
Mean height RMS, mm & 9.6770 & 1.1621 \\
Mean input RMS, m$^2$\,s$^{-2}$ & 0.4582 & 0.4681 \\
\bottomrule
\end{tabular}

\caption{Paired results generated by the supplied code. These are finite-ensemble measurements under the specified disturbance and reference, not universal controller rankings.}\label{tab:uncertainty}
\end{table}
\begin{figure}[tb]
\centering\includegraphics[width=\linewidth]{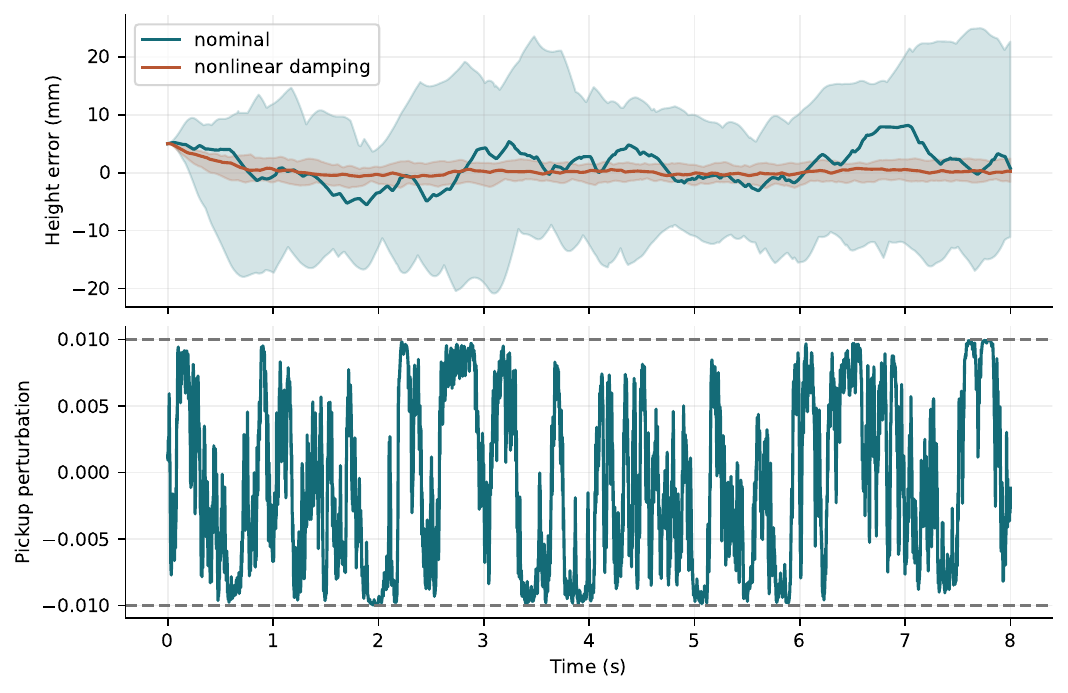}
\caption{Bounded coefficient uncertainty. The curves show empirical medians and pointwise 5th--95th percentile bands across 20 realizations, not confidence intervals for a population mean. The lower panel shows seed 0.}\label{fig:uncertainty}
\end{figure}
The mean RMS height error decreases by \JitterReduction\% in this experiment. The damped trajectories have minimum normalized pickup tension \MinSimTension\,$\mathrm{m^2\,s^{-2}}$, minimum crown speed \MinSimCrown\,$\mathrm{m\,s^{-1}}$, maximum $V$ of \MaxSimV\,$\mathrm{m^2\,s^{-2}}$, and maximum absolute input \MaxSimInput\,$\mathrm{m^2\,s^{-2}}$. These values describe the sampled trajectories. The continuous-time robust guarantees for the ideal model instead follow from Theorem~\ref{thm:robust} and the analytical bounds checked before simulation. The absence of observed failures in a finite set of trials would not by itself establish those guarantees.

\subsection{A baseline and actuator-lag sensitivity}
\begin{figure}[tb]
\centering\includegraphics[width=\linewidth]{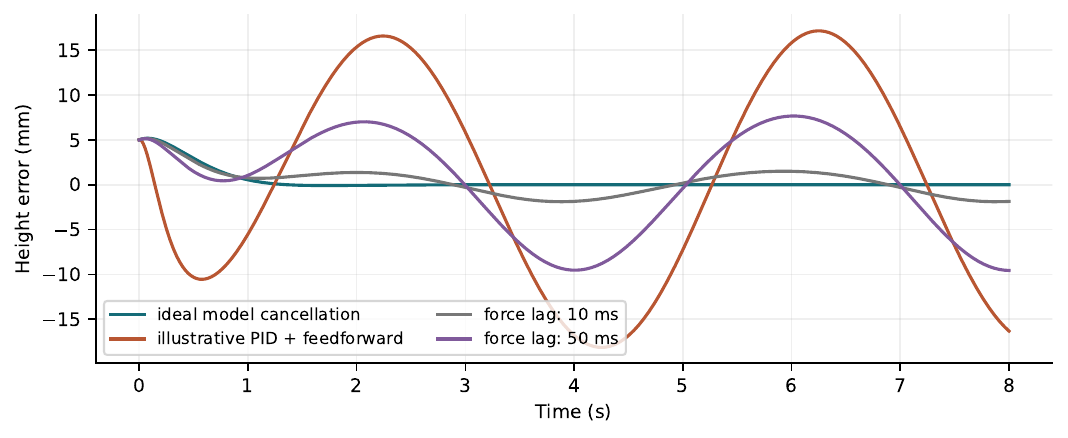}
\caption{Implementation comparisons on the nominal plant. An illustrative PID with steady-state feedforward and two force-lag cases are included. The latter alter the plant and lie outside the proved ideal-servo tracking theorem.}\label{fig:implementation}
\end{figure}

For context, the script also applies the illustrative controller
\begin{equation}
 u=u_r+20(h-r)+8\int_0^t(h-r)\dd s+2(q-\dot r),\label{eq:pid}
\end{equation}
with actuator clipping and integrator clamping. Its gains carry the units implied by the input and signals and are chosen for illustration, without an optimality claim. Separately, the ideal nonlinear command is passed through $\tau\dot u=\un-u$ with $\tau=10$ and $50\,\mathrm{ms}$. Figure~\ref{fig:implementation} displays the resulting errors. These cases illustrate how model cancellation and force bandwidth affect the response. The augmented lag model requires a new invariant-region analysis. Likewise, comparison with a systematically tuned PID or a controller designed for actuator dynamics would require a separate study.

\FloatBarrier
\clearpage\section{Zero gravity: transport obstruction and controlled circulation}\label{sec:zerog}
\subsection{Which experiment are we modeling?}
Setting $g=0$ in a differential equation is meaningful only after specifying the boundary conditions. The vertical fountain has a source, a moving turn, a sink, and constitutive laws for pickup and deposition. A floating closed chain has none of these boundaries and therefore requires a different mechanical model. In a local freely falling frame, its body force can be approximated as zero if tidal forces, residual acceleration, air motion and contact are neglected. This is an idealization of microgravity. It does not require the gravitational field around an orbiting apparatus to vanish.

Our extension therefore addresses two questions. First, how much of the original positive-flow guarantee survives when the body force is removed and the same boundary laws are imposed? Second, what can a speed regulator accomplish for a closed inextensible loop with no body force? For the open system, we find an obstruction that no pickup controller can overcome. For the closed loop, we construct controlled motions within an explicitly described family. Both conclusions concern these particular models, and neither establishes stability for all spatial chain configurations.

The distinction has an established mechanical basis. Dehadrai and Hanna~\cite{dehadrai2026} place the gravity-free lariat, the open chain fountain and the driven string shooter in a common account of axially moving strings. Their analysis shows how gravity, drag, bending stiffness and the drive conditions affect the possible steady shapes. For the closed-loop problem, we take the gravity-free lariat family and ask which changes in circulation speed are compatible with its continuum balance and with prescribed force and tension limits.

\begin{figure}[htbp]\centering
\includegraphics[width=\linewidth]{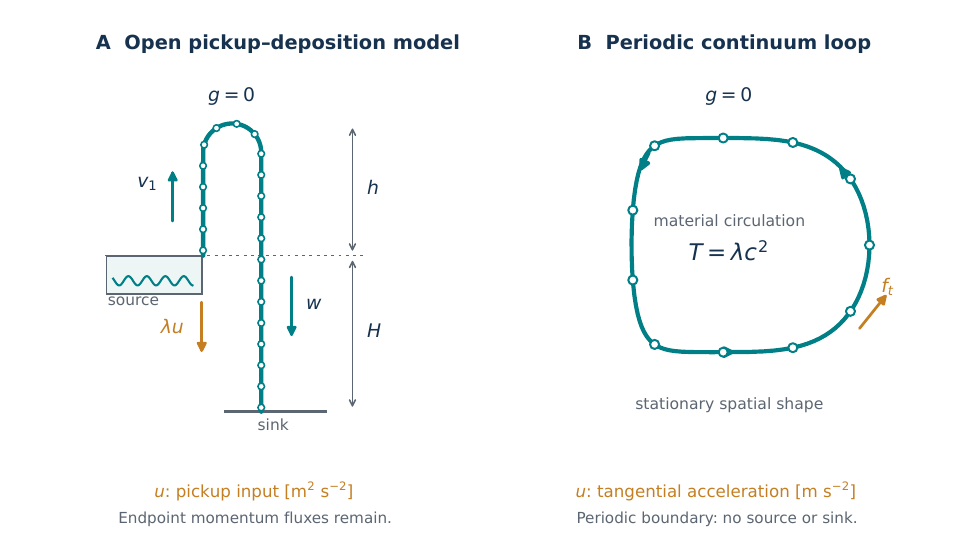}
\caption{Different boundaries require different models. The left panel shows the original open transport idealization with $g=0$. The right panel shows a closed, freely floating continuum loop, with markers indicating material motion along a stationary spatial curve. The input in the right panel is a distributed tangential acceleration, whereas the left panel uses a pickup force.}\label{fig:zerog-models}
\end{figure}

\subsection{An obstruction that the pickup actuator cannot remove}
Consider \eqref{eq:model} with $g=0$, $H>0$, $0\le\beta<1/2$, and $h,v_1,w>0$. The coordinate $h$ and the separation $H$ still refer to the prescribed two-leg geometry, whose direction is now specified independently of gravity. Write
\begin{equation}
 P=(H+h)w,\qquad s=(v_1+w)/2.\label{eq:zg-P}
\end{equation}
Here $\lambda P$ is the magnitude of the falling-leg momentum. Although the pickup input influences the state indirectly, it does not appear explicitly in the derivative of $P$.

\begin{theorem}[Finite duration of a uniform positive crown-speed margin]\label{thm:zg-budget}
Let a classical solution of the zero-gravity open model remain in $\D$ on $[0,t]$. If $s(\tau)\ge s_{\min}>0$ throughout this interval, then
\begin{equation}
 \dot P=-\frac14v_1^2-\left(\frac34-\beta\right)w^2
 \le-\frac14(v_1^2+w^2)\le-\frac12s^2,\qquad
 t<\frac{2(H+h(0))w(0)}{s_{\min}^2}.\label{eq:zg-budget}
\end{equation}
Consequently no solution of this open model can retain both its positive physical domain and a fixed positive crown-speed margin for all forward time, irrespective of the pickup-input law.
\end{theorem}
\begin{proof}
The product rule, $\dot h=(v_1-w)/2$, and the falling-leg equation give
\[
 \dot P=\frac{(v_1-w)w}{2}-\frac{(v_1+w)^2}{4}+\beta w^2
 =-\frac14v_1^2-\left(\frac34-\beta\right)w^2.
\]
Since $\beta<1/2$, the coefficient of $w^2$ is greater than $1/4$. Also $(v_1-w)^2\ge0$ implies $v_1^2+w^2\ge2s^2$. Integrating yields $P(t)\le P(0)-s_{\min}^2t/2$. Positivity of $(H+h(t))w(t)$ makes the time bound strict. The argument applies on every finite interval on which the stated hypotheses hold and thus does not assume global existence.
\end{proof}
The bound limits the time for which all the hypotheses can hold. It does not determine the exact exit time, since a trajectory may lose its prescribed speed margin much earlier. The conclusion depends on this open geometry and its deposition law. An actuator at the other end would change the momentum balance and require a new derivation.

The failure is also visible in the lower speed-face condition of Theorem~\ref{thm:tube}:
\[
 0=g(H+h_{\min})\ge(1-\beta)w_{\min}^2+2Qw_{\min}+Q^2.
\]
For $Q\ge0$ and $w_{\min}>2Q$, the right side is strictly positive. This particular invariant-region construction cannot survive at $g=0$. The momentum argument is stronger: it rules out the desired all-time uniform margin for \emph{any} pickup law within the unchanged open model.

\subsection{What remains possible: constant height with decaying flow}
A stationary height does not require a stationary positive flow. Fix $r>0$ and $w_0>0$, and set
\begin{equation}
 b=\frac{1-\beta}{H+r}>0,\quad
 w(t)=v_1(t)=\frac{w_0}{1+bw_0t},\quad h(t)=r,\quad
 u(t)=(\alpha+rb)w(t)^2.\label{eq:zg-decay}
\end{equation}
\begin{proposition}[Exact constant-height decay]\label{prop:zg-decay}
For $0\le\alpha<1/2$ and $0\le\beta<1/2$, \eqref{eq:zg-decay} solves the zero-gravity open model for every $t\ge0$. Both leg speeds and the pickup tension are positive at every finite time, while $w(t)\to0$. The input is passive braking, and
\begin{equation}
 T_P/\lambda=(1+rb)w^2>0,\qquad
 \dot E/\lambda=-(1/2-\alpha)v_1^3-(1/2-\beta)w^3-u v_1.\label{eq:zg-energy}
\end{equation}
There is no equilibrium with $h,v_1,w>0$ at $g=0$, even with a constant signed pickup input.
\end{proposition}
\begin{proof}
Differentiation gives $\dot w=-bw^2$. With $v_1=w$, the crown speed is $w$, the kinematic equation gives $\dot h=0$, and the falling equation gives $\dot w=-(1-\beta)w^2/(H+r)$. Substituting the displayed input in the rising equation gives the same derivative. The pickup formula gives the stated tension. Equation~\eqref{eq:zg-energy} is the energy identity of Proposition~\ref{prop:energy} at $g=0$. At any equilibrium, the kinematics force $v_1=w$, after which $\dot w=-bw^2<0$ contradicts positive flow.
\end{proof}
The solution exists globally in the ideal model with an inexhaustible source, although its speed and tension have no uniform positive lower bounds. A finite source gives an additional time restriction.
\begin{corollary}[Exact depletion of a finite source]\label{cor:zg-depletion}
For the motion in Proposition~\ref{prop:zg-decay}, let $\ell_0>0$ be the available length and define $\ell(t)=\ell_0-\int_0^t w(\tau)\,\dd\tau$. Then
\begin{equation}
 \ell(t)=\ell_0-b^{-1}\log(1+bw_0t),\qquad
 t_{\rm dep}=\frac{\exp(b\ell_0)-1}{bw_0}>0.\label{eq:zg-depletion}
\end{equation}
For every $t\ge0$, $\ell(t)>0$ if and only if $t<t_{\rm dep}$, and $\ell(t_{\rm dep})=0$. Also $\dot\ell=-w$. Thus the positive-supply part of this exact mechanical motion is precisely $0\le t<t_{\rm dep}$.
\end{corollary}
\begin{proof}
The derivative of $b^{-1}\log(1+bw_0t)$ is $w_0/(1+bw_0t)$. Its value at zero is zero, so the fundamental theorem of calculus gives the integral. Since $b,w_0>0$, strict monotonicity of the logarithm gives the stated equivalence and the unique zero of the remaining supply.
\end{proof}
The model does not describe subsequent motion after the source is exhausted. A similar loss of positive flow occurs in the equilibrium for positive gravity: its speed is proportional to $\sqrt g$ and tends to zero as $g\downarrow0$.

\begin{figure}[htbp]\centering
\includegraphics[width=\linewidth]{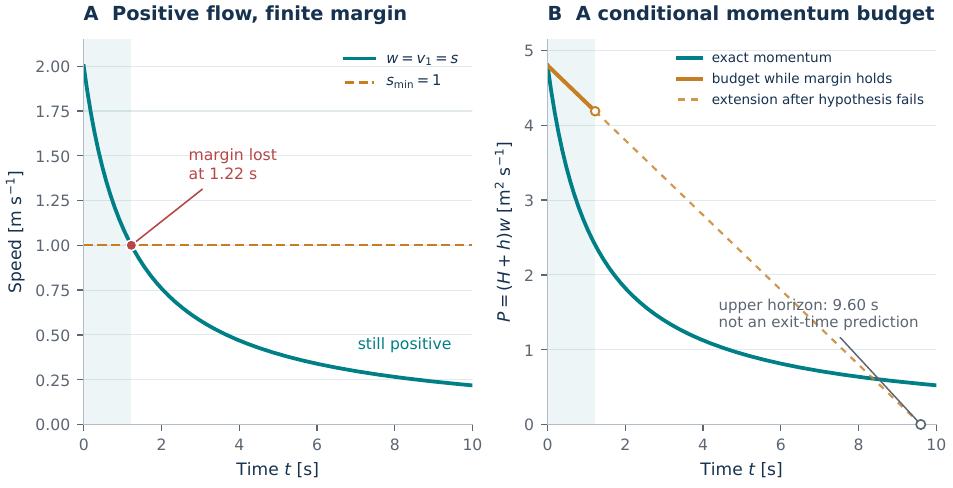}
\caption{Analytical zero-gravity open-model example: $H=2$, $r=0.4$, $w_0=2$, $\alpha=0.12$, $\beta=0.02$ in consistent SI units, with $s_{\min}=1$. The exact constant-height trajectory loses a selected speed margin while remaining positive. The affine momentum budget applies only until that margin is lost. Its zero gives an upper bound on this duration, rather than a simulated exit time.}\label{fig:zg-decay}
\end{figure}

\subsection{The closed-loop continuum model and an exact family}
Let $a\in\mathbb R/L\mathbb Z$ label material arclength, $X(a,t)\in\mathbb R^2$ or $\mathbb R^3$ denote position, and $T(a,t)$ denote tension. For a perfectly flexible inextensible string of density $\lambda>0$, use
\begin{equation}
 \lambda X_{tt}=\partial_a(TX_a)+f,\qquad |X_a|=1,\qquad
 X(a+L,t)=X(a,t).\label{eq:zg-string}
\end{equation}
This is the standard inertial string equation. The unforced motion that motivates our construction is classical: a closed string can circulate at constant speed along a fixed spatial curve, with uniform tension supplying the normal acceleration. Healey and Papadopoulos~\cite{healey1990} give a direct treatment, including arbitrary closed space curves and an extension to nonlinearly elastic strings. Hanna and Santangelo~\cite{hanna2013} explain why a free end changes the problem, since the endpoint condition cannot generally support the positive uniform tension of a circulating loop. Moving pickup or contact discontinuities also require separate jump balances~\cite{hanna2015}.

We use this known family to formulate a controlled motion. The spatial curve will remain fixed, while a distributed tangential force changes the circulation speed. The reconstruction below identifies the force field required for that change. The subsequent tracking theorem combines the resulting scalar equation with positive-speed, tensile-load and actuator bounds.

Choose a $C^2$, $L$-periodic curve $R$ with $|R'(a)|=1$. For a smooth circulation speed $c(t)$ and phase $\theta$ satisfying $\dot\theta=c$, set
\begin{equation}
 X(a,t)=R(a+\theta(t)),\qquad T(a,t)=\lambda c(t)^2.\label{eq:zg-ansatz}
\end{equation}
Individual material points move along the curve while its image in space stays fixed. For a noncircular curve, this motion differs from rigid-body rotation.

\begin{proposition}[Exact reconstruction with tangential forcing]\label{prop:zg-loop}
Let $\gamma\ge0$ and let $u,d$ be spatially uniform scalar accelerations. If
\begin{equation}
 \dot c=u-\gamma c+d,\qquad
 f(a,t)=\lambda\bigl(u(t)-\gamma c(t)+d(t)\bigr)R'(a+\theta(t)),\label{eq:zg-speed}
\end{equation}
then \eqref{eq:zg-ansatz} satisfies \eqref{eq:zg-string}. For continuous $u,d$ this is a classical solution. With locally integrable $u,d$ and locally absolutely continuous $c$, the temporal equation holds almost everywhere. In particular, $u=\gamma=d=0$ gives constant-speed circulation with positive tension when $c\ne0$.
\end{proposition}
\begin{proof}
The chain rule gives $X_a=R'$, $X_{aa}=R''$, $X_t=cR'$, and $X_{tt}=\dot cR'+c^2R''$. Tension is spatially constant, so $\partial_a(TX_a)=\lambda c^2R''$. The remaining tangential equation is exactly \eqref{eq:zg-speed}, and $R$ supplies the unit-speed and periodic constraints. The same differentiation holds almost everywhere for the stated absolutely continuous time dependence. The exceptional set is independent of the material label $a$: outside the null set where the scalar speed equation fails, the spatial chain rule applies at every $a$ simultaneously.
\end{proof}
Normal acceleration $c^2R''$ is supplied by the tension gradient vector, even though the scalar tension is uniform. No gravity is required for this balance. The initial velocity must be tangent with the same material speed everywhere. An arbitrary release, an open end, a localized kick, or nonuniform forcing does not satisfy these hypotheses.

Here $u$ has units of acceleration, whereas the open-model input has units of squared speed. Distributed tangential actuation assumes an ideal mechanism acting along the loop. Localized forcing leads to a different mechanical problem. O'Reilly and Varadi~\cite{oreilly2004} study axial motion of a closed string with a singular supply of momentum, and Hanna~\cite{hanna2015} derives the associated form of jump balances for moving discontinuities. A drive at one contact point must satisfy such local balances and cannot simply be identified with the distributed input $u$ used here.

\subsection{What the regulator regulates}
The controlled output is now circulation speed $c$, not apex height. Given a $C^1$ reference $r(t)$ and gain $k>0$, use
\begin{equation}
 u=\dot r+\gamma r-k(c-r),\quad e=c-r,\quad \mu=\gamma+k>0.
 \qquad \dot e=-\mu e+d.\label{eq:zg-control}
\end{equation}
The reference derivative supplies acceleration, $\gamma r$ compensates drag at the reference, and $-ke$ corrects speed error. Passive drag contributes to the error-decay rate: it is $\mu=\gamma+k$, not merely $k$. This regulator does not select a preferred geometry from the family of possible curves.

The scalar feedback calculation uses the standard integrating-factor method for a stable first-order equation. Its role here is to produce a speed history that can be reconstructed as a continuum motion and checked against physical limits. The theorem makes this connection explicit for measurable disturbances, including existence and uniqueness of the scalar solution and simultaneous bounds on speed, tension and command. The classical circulation family and the elementary scalar decay estimate are the starting points of that argument.

\begin{theorem}[Global speed tracking and physical margins on the exact family]\label{thm:zg-tracking}
Let $\lambda>0$, $\gamma\ge0$, $k>0$, $A,D\ge0$, $c_0\in\mathbb R$, and $r\in C^1([0,\infty))$. Suppose $r_{\min}\le r(t)\le r_{\max}$, $|\dot r(t)|\le A$, and $d$ is measurable with $|d(t)|\le D$ almost everywhere, for $t\ge0$. Let
\begin{equation}
 \rho=\max\{|c_0-r(0)|,D/\mu\},\qquad r_{\min}>\rho.\label{eq:zg-rho}
\end{equation}
The scalar closed-loop equation has a global locally absolutely continuous solution with $c(0)=c_0$, unique among all such solutions on $t\ge0$ with the same initial value, and
\begin{align}
 |e(t)|&\le |e(0)|e^{-\mu t}+\frac D\mu(1-e^{-\mu t})\le\rho,\label{eq:zg-bound}\\
 0<r_{\min}-\rho&\le c(t)\le r_{\max}+\rho,\label{eq:zg-speedbounds}\\
 \lambda(r_{\min}-\rho)^2&\le T(t)\le\lambda(r_{\max}+\rho)^2,\qquad
 |u(t)|\le A+\gamma r_{\max}+k\rho.\label{eq:zg-forcebounds}
\end{align}
Consequently an acceleration limit $U_{\max}$ and a tensile load limit $T_{\max}$ are respected whenever the respective right sides in \eqref{eq:zg-forcebounds} do not exceed these limits. Together with \eqref{eq:zg-ansatz}, this constructs a global solution within the exact circulation family. It is not a global existence theorem for arbitrary initial data of \eqref{eq:zg-string}.
\end{theorem}
\begin{proof}
Variation of constants constructs
\[
 e(t)=e(0)e^{-\mu t}+\int_0^t e^{-\mu(t-\tau)}d(\tau)\,\dd\tau.
\]
On each compact interval the integral exists and is absolutely continuous, and differentiation gives the equation almost everywhere. The difference of any two solutions with the same initial value solves $\dot z=-\mu z$, so uniqueness follows by multiplication by $e^{\mu t}$. Taking absolute values and integrating the positive kernel gives \eqref{eq:zg-bound}. Since its right side is a convex combination of $|e(0)|$ and $D/\mu$, it is at most $\rho$. Adding the reference bounds proves \eqref{eq:zg-speedbounds}. These bounds are positive, so they can be squared to obtain the tension bounds. The triangle inequality applied to the control law gives the command bound. Extend $c$ continuously to negative times by its initial value and set $\theta(t)=\theta_0+\int_0^t c(s)\,\dd s$ using that extension. The ordinary fundamental theorem of calculus gives $\dot\theta=c$ at every forward time, including zero. Proposition~\ref{prop:zg-loop} then reconstructs $X$, with its equation valid on one full-measure time set for all material labels. The choice of negative-time extension does not restrict the forward solution or its uniqueness.
\end{proof}
For $D=0$ the error decays exponentially, while persistent bounded disturbances give an ultimate bound $D/\mu$. Increasing $k$ improves that bound but can increase the input required initially. Saturation changes \eqref{eq:zg-control}, so the tracking estimate applies only when the sufficient input limit holds or a separate saturation analysis justifies it. Brownian forcing falls outside the bounded measurable disturbance assumption. The earlier SDE results for the open fountain therefore do not establish a stochastic theorem for this closed-loop PDE.

\begin{figure}[htbp]\centering
\includegraphics[width=\linewidth]{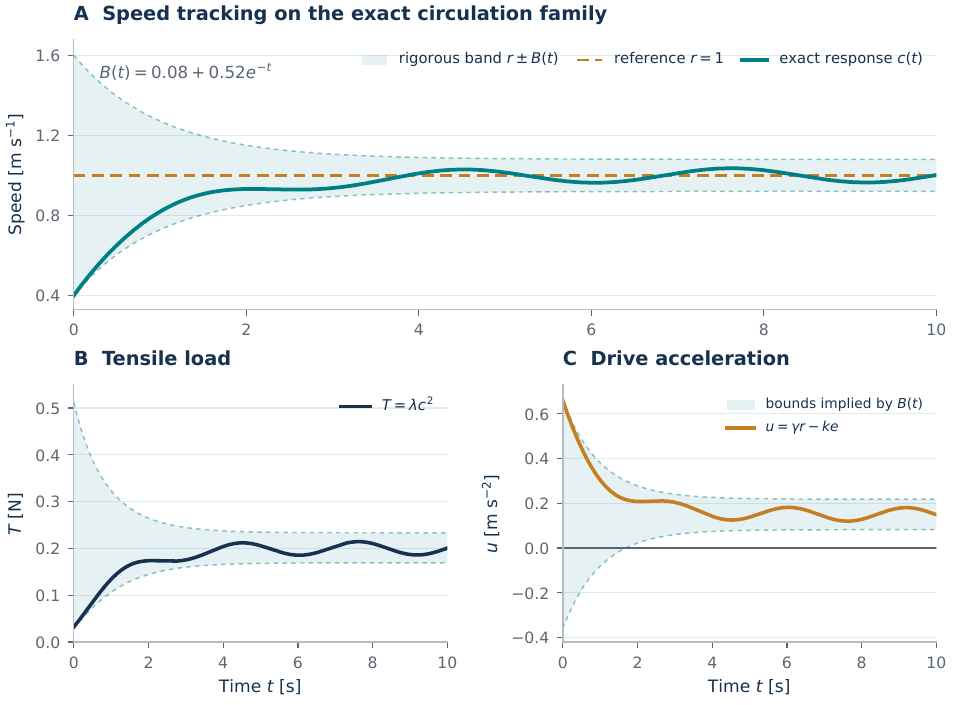}
\caption{An explicit speed-regulation example on the exact family: $r=1$, $c(0)=0.4$, $\gamma=0.15$, $k=0.85$, $d=0.08\sin(2t)$, $\lambda=0.2$. The analytical solution stays inside the rigorous error envelope, with positive bounded tension and drive. The geometry is arbitrary within the smooth unit-speed periodic family and is not stabilized by this scalar law.}\label{fig:zg-tracking}
\end{figure}

\subsection{Energy, reaction torque, and the obstruction to shape attraction}
The kinetic energy on the family is
\begin{equation}
 E=\frac12\lambda Lc^2,\qquad
 \dot E=\lambda Lc(u-\gamma c+d).\label{eq:zg-loop-energy}
\end{equation}
For measurable forcing, the derivative identity holds almost everywhere. Internal tension does no net work on this closed inextensible motion: differentiating $|R^\prime|^2=1$ gives $R^\prime\!\cdot R^{\prime\prime}=0$, so $\partial_a(TX_a)\cdot X_t=0$ pointwise. The spatial integral of external power is $\int_0^L f\cdot X_t\,\dd a=\lambda Lc(u-\gamma c+d)$, which agrees with the energy derivative. Drag removes energy, while the external drive can supply or extract it. Because $\int_0^L R'\,\dd a=0$, a spatially uniform tangential drive has zero resultant force, although its torque may be nonzero. For a planar loop with signed area $\mathcal A:=\tfrac12\int_0^L R\times R^\prime\,\dd a$ (interpreting the cross product as its scalar normal component),
\[
 J=2\lambda\mathcal A c,\qquad
 \dot J=2\lambda\mathcal A(u-\gamma c+d).
\]
These formulas use the physical definitions $J(t)=\int_0^L\lambda X(a,t)\times X_t(a,t)\,\dd a$ and $\tau(t)=\int_0^L X(a,t)\times f(a,t)\,\dd a$. Periodicity makes the spatial integral invariant under the material phase shift. Substituting $X_t=cR'$ and the tangential force therefore gives $J=2\lambda\mathcal A c$ and $\tau=2\lambda\mathcal A(u-\gamma c+d)$. The scalar speed equation gives $\dot J=\tau$, classically or almost everywhere according to the regularity of $c$. The zero resultant force follows separately by integrating $R'$ over a period. When $\mathcal A\ne0$, accelerating circulation requires an angular-momentum exchange with an actuator or environment. The supporting device and its attitude dynamics remain outside the current model.

\begin{proposition}[Geometric obstruction to shape attraction]\label{prop:zg-shape-obstruction}
Work in a fixed spatial frame and identify the shape with its geometric image, independently of its material parameterization. For a continuous $L$-periodic curve $R$ with $L>0$, let $K(t)=\{R(a+\theta(t)):a\in\mathbb R\}$. For any phase history, $K(t)=K(0)$. Hence, for every nonempty compact target $K_*$, the Hausdorff distance satisfies
\[
 d_H(K(t),K_*)=d_H(K(0),K_*).
\]
It tends to zero as $t\to\infty$ if and only if $K(0)=K_*$. In particular, two admissible curves of the same length with distinct geometric images cannot both be attracted to one fixed target by the scalar circulation law. With the same initial speed, reference and disturbance, their speed histories agree.
\end{proposition}
\begin{proof}
Translation of the material variable is a bijection of $\mathbb R$, so the range of $a\mapsto R(a+\theta(t))$ is the range of $R$ at every time. Continuity and periodicity make this range nonempty and compact. Hausdorff distance between two such sets is zero exactly when the sets agree; a constant positive distance cannot converge to zero. Equality of the speed histories follows from uniqueness of the scalar closed-loop equation.
\end{proof}
This result concerns attraction to a spatial image in a fixed frame; it does not quotient shapes by rigid motions. Stability to general transverse perturbations remains a separate question that requires the full constrained spatial dynamics.

There are substantial stability results for related loop models. Healey~\cite{healey1996} combines conserved energy and circulation in a Lyapunov function to establish orbital stability of axial motions of nonlinearly elastic loops under conditions on the material and axial speed. The comparison is with a family of relative equilibria, so that conclusion is compatible with the absence of attraction to one prescribed shape. For a driven string shooter with bending stiffness, Kocbay~\cite{kocbay2026} studies the spectrum of linearized spatial perturbations in a Kirchhoff-rod model. These works address aspects of spatial stability that our scalar error estimate leaves open. Applying their conclusions to the present forced inextensible loop would require agreement of the constitutive law, forcing, boundary conditions and perturbation class.

\subsection{Simulation protocol and evidential scope}\label{sec:zg-simulation}
The companion presents three mathematical levels. The original Earth mode integrates the established three-state model. The open zero-gravity mode uses the same geometry and lets the reader choose positive initial speeds. Substituting $g=0$ in the positive-gravity equilibrium formula would instead initialize the system on the zero-flow boundary. The momentum-budget diagnostic tests Theorem~\ref{thm:zg-budget} only while its hypotheses hold.

The exact-loop view renders material tracers along a periodic curve and integrates the scalar circulation equation. These tracers sample a continuum solution. Equal material spacing measures arclength and does not imply fixed finite chord lengths on a noncircular curve. Reference changes, drag and bounded homogeneous disturbances illustrate Theorem~\ref{thm:zg-tracking}. The rendered geometry stays fixed by the model's ansatz, so this view cannot be used as evidence of shape stability.

A separate elastic bead-ring experiment evolves Cartesian bead positions with tension-only links. A local kick can deform the ring, and no projection returns it to the reference curve. This regularized finite-dimensional model permits extension and slack, which places it outside the inextensible continuum PDE and its theorem. Momentum, energy, angular momentum, slack and timestep comparisons serve as numerical diagnostics. On a deformed discrete ring the vertex tangent forces can have a nonzero resultant. The displayed force and momentum diagnostics make this difference from the exactly force-balanced continuum family visible. Because the experiment omits self-contact and bending physics, a self-crossing or highly folded configuration cannot be interpreted as a faithful prediction for a real beaded chain. The formal correspondence for the new results is recorded in Appendix~\ref{app:zg-math}.

\FloatBarrier\clearpage\section{Scope and next developments}
The control interpretation is useful because it identifies both a manipulable boundary quantity and the conditions under which its use preserves the underlying mechanics. A stable height transfer function explains why modest feedback is plausible. The energy balance distinguishes passive suppression from active elevation. The invariant-region inequalities then connect tracking speed to positive transport and force capacity. These are results for the adopted reduced model, with a specified force input and measurable states. The big picture here is moving towards stabilization in 3d that would be useful for real-world applications, related to hovering aicrafts.

The next physical step is to identify the boundary law in an actuated apparatus. Measurements should determine whether the traction device changes the effective pickup coefficient, how much force it can supply at the relevant transport speeds, and which bandwidth is available. Independent discrete-chain calculations or experiments should then test the predictions when finite links, changing packing, and transverse motion are present. Reusing the same ODE for both controller design and simulation is a reproducibility exercise, not such an independent test.

There are also mathematical extensions with clear starting points. Replacing the conservative rectangular bounds by a larger viability region would also allow a constrained reference governor to use that region online. Actuator dynamics and observation errors can be included as additional states or bounded perturbations, with new invariant-set estimates. For stochastic forcing, a physically justified noise model should be selected before choosing almost-sure, moment, or chance-constrained objectives. Each extension would require additional hypotheses and proofs beyond the present theorems. We make no categorical claim that this is the first control-theoretic treatment of the phenomenon.

\clearpage\appendix
\section{Mathematical background and a guide to the proofs}\label{app:math}
This appendix explains the analytical tools used in the proofs, starting from ordinary calculus, finite-dimensional linear algebra and elementary probability. It assumes no prior course in nonlinear or stochastic control. We prove the elementary estimates used in the paper and state the deeper foundational results together with their hypotheses. The new examples and background arguments serve this explanatory purpose. Their correspondence with the formal development is recorded separately in Appendix~\ref{app:formal}, and they do not constitute additional Lean declarations.

\subsection{Functions, trajectories, and the meaning of ``almost everywhere''}\label{app:ac}
A property holds almost everywhere in time if the exceptional set has Lebesgue measure zero. An integral is unchanged when its integrand is modified on such a set. This is the appropriate language for a bounded measurable disturbance, which may jump and need not have a classical derivative.

A scalar function $x:[a,b]\to\R$ is absolutely continuous if, for every $\varepsilon>0$, there is $\delta>0$ such that every finite disjoint family of intervals with total length below $\delta$ has total endpoint variation below $\varepsilon$. Lipschitz functions are absolutely continuous, but the converse need not hold. Local absolute continuity on $[0,\infty)$ means absolute continuity on every compact time interval.

The fundamental theorem of Lebesgue calculus says
\[
x\text{ is AC on }[a,b]
\quad\Longleftrightarrow\quad
x(t)=x(a)+\int_a^t g(s)\,\dd s\text{ for some }g\in L^1(a,b).
\]
In that case $x'$ exists almost everywhere, belongs to $L^1$, and equals $g$ almost everywhere. The vector version follows coordinate by coordinate. The direction from an integral to its derivative uses Lebesgue differentiation. The reverse direction is precisely where absolute continuity is stronger than continuity and almost-everywhere differentiability.

This matters for the mechanical equations. A Carath\'eodory solution means a locally absolutely continuous curve satisfying $\dot x=F(t,x)$ almost everywhere. It is equivalent to
\[
x(t)=x(0)+\int_0^tF(s,x(s))\,\dd s
\]
when the integrand is locally integrable. On a compact state box, the bounded disturbance and reference data provide a uniform bound on $F$, so the model itself supplies the required integrability. Absolute continuity is essential here, since continuity alone does not give a curve the derivative required by the ODE.

For a concrete example, let $d(t)$ switch between two allowed coefficient values at time $t_0$. The integral equation describes a continuous state across $t_0$, even though its acceleration may jump. Integrating the energy estimate therefore does not require the differential equation to hold at that one switching time. The same reasoning treats bounded measurable switching without selecting individual switching times.

If $\Phi$ is continuously differentiable near the compact range of an AC curve, then $\Phi\circ x$ is AC and
\[
\frac{\dd}{\dd t}\Phi(x(t))=D\Phi(x(t))\dot x(t)
\quad\text{almost everywhere}.
\]
The time-dependent version includes the partial derivative in time. This is the chain rule used for $E$, $V$, and the transformed errors. It remains valid when $d$ has no derivative. The MIACE companion~\cite[Sections 2.5 and 7.1]{tomov13} provides a complementary treatment of the same real-analysis foundation.

\subsection{Integral inequalities and why exponential factors work}\label{app:comparison}
Suppose an AC nonnegative function satisfies $\dot V\le-cV+b(t)$ almost everywhere, where $c>0$ and $b$ is locally integrable. The product rule gives
\[
\frac{\dd}{\dd t}\{e^{ct}V(t)\}\le e^{ct}b(t)
\]
almost everywhere. Integration, followed by division by $e^{ct}$, yields
\begin{equation}
V(t)\le e^{-ct}V(0)+\int_0^te^{-c(t-s)}b(s)\,\dd s.\label{eq:background-comparison}
\end{equation}
For a constant bound $b(s)\le B$, the integral is $B(1-e^{-ct})/c$. Set $c=2\kappa$ and $B=\bar d^2/(8\kappa_d)$ to recover the robust error envelope. For $b=0$ one obtains nominal exponential decay. If $V(0)\le\Vbar$ and $B/c\le\Vbar$, the right side is a convex combination of two numbers at most $\Vbar$, which proves invariance of the energy cap.

The related Gronwall inequality supplies uniqueness. If a continuous nonnegative function $y$ satisfies $y(t)\le A+L\int_0^ty(s)\dd s$, define $Y(t)=A+L\int_0^ty(s)\dd s$. Then $Y'\le LY$ and $Y(0)=A$, so $y(t)\le Ae^{Lt}$. In particular, $A=0$ forces $y=0$. Applied to the distance between two integral solutions with the same initial value and a state-Lipschitz field, this proves equality on their common interval.

For two exponential rates $a,b>0$, variation of constants produces the convolution
\[
\int_0^te^{-a(t-s)}e^{-bs}\,\dd s
=\begin{cases}(e^{-bt}-e^{-at})/(a-b),&a\ne b,\\te^{-at},&a=b.
\end{cases}
\]
Even in the equal-rate case, $te^{-at}$ is bounded by a constant times $e^{-ct}$ for every $0<c<a$. The polynomial prefactor in the convolution thus explains why the full-state convergence corollary permits a slightly smaller decay rate when the two rates coincide.

\subsection{Local construction, compact continuation, and uniqueness}\label{app:existence}
The measurable-time local theorem needs three ingredients: $F(t,x)$ is measurable in $t$, locally Lipschitz in $x$ uniformly on the relevant compact parameter ranges, and locally bounded by an integrable function. In the present model, the only state singularities occur at the height denominators. On a compact box strictly inside $h>0$, smooth dependence on state and bounded reference/disturbance data provide these conditions.

For a bounded globally Lipschitz field, start with the Picard integral operator
\[
(\mathcal T x)(t)=x_0+\int_0^tF(s,x(s))\,\dd s.
\]
On a short interval, $\|\mathcal T x-\mathcal T y\|_\infty\le LT\|x-y\|_\infty$, so it is a contraction when $LT<1$. Measurability suffices for the integral. Iteration converges in the complete space of continuous curves, and the limit solves the integral equation. One can alternatively estimate the $n$th iterate by $(LT)^n/n!$, avoiding a restriction on the final finite horizon. Local state restrictions are handled by a bounded Lipschitz extension or by choosing a ball that Picard iteration does not leave.

To continue a locally constructed solution, a uniform bound $\|\dot x\|\le M$ gives $\|x(t)-x(s)\|\le M|t-s|$. If a finite endpoint $\tau$ existed while $x$ remained in a fixed compact subset of the regular domain, the curve would have a limit $x_*$. The local theorem at $(\tau,x_*)$ extends it. The two pieces share their endpoint and their integral equation, so the concatenation is again AC. This contradicts maximality.

In the paper's global theorem, estimates on each existing finite prefix justify continuation to later times. This order closes the existence argument, whereas boundedness of an assumed all-time trajectory would leave its construction unresolved. The Lean proof first projects onto an outer state box to obtain a bounded extension, then proves that the solution stays in an inner box where the projection is inactive. The projection therefore serves the existence proof without saturating the commanded force.

To establish uniqueness among unrestricted competitors, observe that continuity keeps a competitor in the regular domain initially. On every common prefix, the same error and face estimates place it in the invariant box, preventing it from reaching a singular boundary where the local uniqueness argument would cease to apply. Thus the formal theorem derives the competitor's box membership rather than requiring it as a hypothesis.

\subsection{Invariant intervals and weak inward inequalities}\label{app:faces}
For a scalar ODE, the intuitive condition at a lower boundary $b$ is $f(t,b)\ge0$. Equality is allowed: the trajectory may touch or slide along the boundary. A strict first-derivative contradiction is therefore not enough. Here the special quadratic structure gives an exact comparison.

For the falling-speed equation, with known continuous $h,q$, set
\[
f(t,w)=g-\frac{(1-\beta)w^2+2q(t)w+q(t)^2}{H+h(t)}.
\]
Subtract its value at a constant $b$:
\[
f(t,w)-f(t,b)=a_b(t)(w-b),\quad
a_b(t)=-\frac{(1-\beta)(w+b)+2q}{H+h}.
\]
At the lower boundary define $y=b-w$. Then $\dot y=a_b y-f(t,b)\le a_by$. Multiplying by $\exp(-\int_0^ta_b)$ shows $y(t)\le0$ whenever $y(0)\le0$. The upper-boundary proof uses $y=w-b$ and $f(t,b)\le0$. Integrability of $a_b$ on finite intervals suffices, regardless of its sign. The comparison therefore establishes invariance under the stated weak boundary inequalities.

For a general multidimensional domain, an inward condition has to be checked at every relevant boundary and must be paired with an appropriate existence theorem. A crown constraint is only one boundary function. The counterexample at $w=0$ demonstrates why protecting its sign does not automatically protect the entire physical domain.

\subsection{The cubic stability test and nonlinear local decay}\label{app:hurwitz}
Consider $p(s)=s^3+a_2s^2+a_1s+a_0$. For positive coefficients, the cubic Hurwitz criterion is $a_2a_1>a_0$. Its Routh array has first column
\[
1,\quad a_2,\quad \frac{a_2a_1-a_0}{a_2},\quad a_0.
\]
All entries are positive exactly under the stated conditions. There is also a useful direct exclusion proof. If a root is real and nonnegative, all terms in $p$ are nonnegative and $a_0>0$, a contradiction. For a nonreal root $s=x+iy$ with $x\ge0$, divide the imaginary part by $y$ to obtain
\[
y^2=3x^2+2a_2x+a_1.
\]
Substituting into the real part gives
\[
0=-8x^3-8a_2x^2-2(a_1+a_2^2)x+(a_0-a_2a_1)<0,
\]
another contradiction. Thus no root lies in the closed right half-plane. The paper's factorization supplies all the needed signs from its physical parameter range.

The standard Lyapunov argument for passing from a stable linearization to nonlinear local stability is developed in \cite[Chapter 4]{khalil}. Write the perturbation equation as $\dot y=Jy+R(y)$, with $\|R(y)\|/\|y\|\to0$. A Hurwitz matrix has an exponentially decaying matrix exponential. For any positive definite $Q$, define
\[
P=\int_0^\infty e^{J^\top t}Qe^{Jt}\,\dd t.
\]
The integral converges, $P$ is positive definite, and differentiation under the integral gives $J^\top P+PJ=-Q$. The quadratic Lyapunov function $W=y^\top Py$ satisfies
\[
\dot W=-y^\top Qy+2y^\top PR(y).
\]
Choose a small ball where $\|R(y)\|\le\varepsilon\|y\|$ and $2\|P\|\varepsilon<\lambda_{\min}(Q)$. Then $\dot W\le-c\|y\|^2\le-cW/\lambda_{\max}(P)$ for a positive $c$. A sufficiently small sublevel of $W$ lies inside that ball and is invariant. Comparison and continuation give local exponential stability.

The matrix integral gives a standard explanation of the nonlinear stability argument. In the checked development, a positive definite quadratic Lyapunov function is constructed for the cubic companion representation and transferred to the mechanical Jacobian before the nonlinear remainder is estimated. Both presentations establish the same statement. The correspondence table records this agreement of conclusions without suggesting that each displayed step reproduces the formal proof literally.

\subsection{Relative degree, zero dynamics, and stable inversion language}\label{app:zeros}
For a linearized realization $(J,B,C)$ with zero initial perturbation, the transfer function is $C(pI-J)^{-1}B$. It summarizes the input/output relation, while internal stability concerns the full state realization. Algebraic cancellation in the transfer fraction can hide a state, so the noncancellation calculation is needed here to rule out cancellation of the one finite zero with a pole.

The relative-degree and zero-dynamics concepts used here follow the usual nonlinear input/output formulation~\cite[Chapter 13]{khalil}. Relative degree two means that the input appears in the second derivative of height, not its first. The nonlinear identity $\dot q=A_0$ is therefore the natural point at which to assign a desired output acceleration. But only two of the three states are fixed by height and its derivative. Setting their errors to zero leaves the zero dynamics in $w$.

For constant $r$, define $c=(1-\beta)/(H+r)>0$. Then
\[
\dot w=-c(w+w_r)(w-w_r).
\]
At $w=0$, the derivative is $g>0$, so the positive branch is invariant. For $\eta=w-w_r$, the exact identity is $\dot\eta=-c(w+w_r)\eta$, and $w+w_r\ge w_r$. Hence
\[
|w(t)-w_r|\le |w(0)-w_r|e^{-cw_rt}.
\]
This proves attraction on that branch. Allowing a nonzero $q$ introduces extra terms in the boundary derivative, and boundedness of $q$ alone need not preserve positivity. The required quantitative bound must also account for $h$, which is why the zero-dynamics result does not amount to an unrestricted input-to-state stability assertion.

\subsection{Constraint geometry, projection, and reference shaping}\label{app:constraints}
For a closed nonempty interval $[L,U]$, minimizing $(u-u_{\rm des})^2$ has the unique solution $u^*=\min\{U,\max\{L,u_{\rm des}\}\}$. If $u_{\rm des}<L$, distance increases as one moves right from $L$, making the left endpoint optimal. The case $u_{\rm des}>U$ is symmetric, while inside the interval the minimum is zero. When $L>U$, the feasible set is empty and the formula no longer solves the stated optimization problem.

This elementary convex geometry explains the robust crown filter. The full reference/gain/state design problem remains more complicated because the sufficient force bounds contain products, roots and inverse heights. The public 2DOF guide~\cite[Sections 2--4]{tomov14} gives related examples in which affine dependence in two weights produces convex sets. Here we use the same practice of writing the feasible set explicitly, without asserting convexity of the full design problem.

For the quintic reference profile $p(\sigma)=10\sigma^3-15\sigma^4+6\sigma^5$, its first derivative is $p'=30\sigma^2(1-\sigma)^2$, maximized at $1/2$ with value $15/8$. Its second derivative is $p''=60\sigma(1-\sigma)(1-2\sigma)$. The extrema of $p''$ occur at $\sigma=(3\pm\sqrt3)/6$ and have magnitude $10/\sqrt3$. Multiplication by the height increment and the time scalings $T^{-1},T^{-2}$ proves~\eqref{eq:quintic}. Both derivatives vanish at the appropriate endpoints for a $C^2$ constant extension. A longer transition reduces derivative demand, while endpoint force requirements remain unchanged.

\subsection{Probability spaces, adaptedness, and the It\^o integral}\label{app:ito}
A filtered probability space $(\Omega,\mathcal F,(\mathcal F_t),\mathbb P)$ records the information available up to time $t$. A process is adapted if its value at $t$ is $\mathcal F_t$-measurable. Brownian motion has continuous paths, independent Gaussian increments relative to this information, and $W_t-W_s$ has variance $t-s$. A continuous adapted process is predictable, which makes its present coefficient admissible in a nonanticipating stochastic integral.

For a simple predictable coefficient $c_s=\sum_jc_j\mathbf1_{(t_j,t_{j+1}]}(s)$, with $c_j$ known at $t_j$, define
\[
\int_0^T c_s\,\dd W_s=\sum_jc_j(W_{t_{j+1}}-W_{t_j}).
\]
Independence and centering remove mixed terms in its second moment. Thus
\[
\mathbb E\left|\int_0^Tc_s\,\dd W_s\right|^2
=\mathbb E\int_0^T|c_s|^2\,\dd s.
\]
Completion in this norm defines the integral for predictable square-integrable coefficients. The resulting integral is a martingale, whose conditional future expectation equals its current value. For locally square-integrable coefficients, stopping where the square integral is bounded gives a local martingale. The bounded stopped integrands in this paper provide the additional integrability needed for zero expectation, a property that an arbitrary local martingale need not have.

Equality at each fixed time almost surely is weaker than equality for all times on one event of probability one. For continuous versions, equality on the countable set of rational times extends to all times. This distinction is crucial when paths are compared up to a random exit, or when a small-time limit quantifies over every sufficiently small real time. The formal construction explicitly proves the simultaneous identities it uses.

\subsection{It\^o's formula, the generator, and stopping}\label{app:generator}
For a continuous diffusion $\dd X=b\dd t+\Sigma\dd W$ and a twice continuously differentiable function $\Phi$, It\^o's formula is
\[
\dd\Phi(X)=\left(D\Phi\,b+\tfrac12\operatorname{tr}(\Sigma\Sigma^\top D^2\Phi)\right)\dd t
+D\Phi\,\Sigma\dd W.
\]
The extra half-Hessian term is forced by Brownian quadratic variation: sums of squared increments converge to elapsed time. Ordinary differentiable paths have no such term. In the error variables here, only $z$ carries diffusion, so the calculation reduces to the one-dimensional second derivative of $z^2/8$. This gives the exact positive correction displayed in~\eqref{eq:generator}.

A stopping time $\tau$ is a random time for which $\{\tau\le t\}$ is known at time $t$. The first exit of a continuous adapted path from an open set has this property under the filtration conventions used in the construction. Continuity gives an almost surely positive exit time when the initial point lies strictly inside, although the positive interval may depend on the path. It need not have a deterministic lower bound shared by all paths.

Stopping keeps the coefficients in a compact region where their bounds apply and makes the martingale terms integrable on a finite horizon. An estimate for $X_{t\wedge\tau}$ therefore concerns the stopped path. Extending it to $X_t$ requires further information about $\tau$. Likewise, a compact exhaustion gives successive exits $\tau_n$ and a maximal lifetime, but passage to the limit $n\to\infty$ requires control of both the exits and integrability.

For bounded Lipschitz coefficients, the SDE Picard iteration replaces one ordinary integral by a stochastic one. It\^o isometry and a martingale maximal inequality control its mean-square supremum. Iteration on short intervals, followed by gluing, gives a continuous strong solution, and analogous difference estimates prove uniqueness in the stated class. Globally Lipschitz extensions on larger compact sets, together with compatibility of the resulting solutions, give the local mechanical construction. Standard treatments such as~\cite{baudoin} develop these stochastic existence tools, while the companion's pinned library development checks the particular route used here.

\subsection{Exponential bounds, Borel--Cantelli, and the freezing argument}\label{app:freezing}
For a bounded predictable $c$, let $M_t=\int_0^tc_s\dd W_s$ and $\langle M\rangle_t=\int_0^tc_s^2\dd s$. The Dol\'eans exponential
\[
Z_t=\exp\{\eta M_t-\eta^2\langle M\rangle_t/2\}
\]
is a mean-one martingale on finite horizons. Boundedness provides the needed exponential integrability. If $|c|\le C$, a maximal estimate applied with $\langle M\rangle_T\le C^2T$ gives
\[
\mathbb P(\sup_{s\le T}M_s\ge x)\le
\exp(-\eta x+\eta^2C^2T/2).
\]
Optimize at $\eta=x/(C^2T)$ and repeat for $-M$ to obtain the two-sided bound used in Theorem~\ref{thm:freezing}. This bound has the correct Gaussian scale even though the integrand is random and state dependent.

The first Borel--Cantelli lemma says that if $\sum_n\mathbb P(A_n)<\infty$, then almost surely only finitely many $A_n$ occur. No independence is needed. For $T_n=\vartheta^n$ and $x_n=C_0C\sqrt{T_n\log n}$, the maximal probabilities are at most $2n^{-C_0^2/2}$, summable when $C_0>\sqrt2$. Since $\log\log(1/T_n)=\log(n\log(1/\vartheta))$, this is precisely the iterated-logarithm scale up to a fixed factor.

The freezing proof uses this coarse bound only for the residual coefficient. Continuity of $a$ gives $|a_s-a_0|\le\varepsilon$ on an initial interval that depends on the path, rather than on a common deterministic interval. Clipping to $[-\varepsilon,\varepsilon]$ supplies a globally bounded coefficient, and stochastic-integral locality identifies its integral with that of the actual residual wherever the coefficients agree initially. Using a countable sequence of clipping levels puts all identities and bounds on one conull event. Sending the level to zero then proves the little-$o$ remainder and completes the passage from the Brownian theorem to the state-dependent mechanical coefficient.

\subsection{The law of the iterated logarithm and its mechanical use}\label{app:lil}
The sharp Brownian LIL at zero is
\[
\limsup_{t\downarrow0}\frac{W_t}{\sqrt{2t\log\log(1/t)}}=1,
\qquad
\liminf_{t\downarrow0}\frac{W_t}{\sqrt{2t\log\log(1/t)}}=-1
\quad\text{almost surely}.
\]
It is a statement about repeated path excursions, stronger than a variance calculation at a fixed time. The upper bound follows from exponential maximal estimates and geometric blocking with parameters approaching their limiting values. The lower bound uses independent increments, Gaussian lower tails and the second Borel--Cantelli lemma to force infinitely many large excursions. The correct Gaussian normalization is necessary to recover the sharp constant. Brownian time inversion transfers the corresponding large-time result to zero. The formal companion proves these ingredients rather than importing an unfinished LIL assertion.

The stochastic-integral freezing theorem transfers both signs to $a_0W_t$ plus a negligible remainder. The deterministic drift is negligible because
\[
\frac{t}{\ell(t)}=\sqrt{\frac{t}{2\log\log(1/t)}}\longrightarrow0.
\]
If one changes the time unit to $t_0$, then $\log\log(t_0/t)/\log\log(1/t)\to1$, so the normalized limit is unchanged. This resolves the dimensional notation inside the logarithm.

At the crown boundary the initial value is zero, and the normal diffusion $a_0=\sigma v_{1,0}^2/(2h_0)$ is nonzero. Therefore both signs occur before every positive time, including before the genuine local exit. This proves non-invariance of the closed boundary constraint under the white-noise model. Interior exit probabilities require a separate analysis, as does the physical adequacy of white noise for any particular disturbance.

\subsection{A compact glossary of logical distinctions}
\begin{longtable}{>{\raggedright\arraybackslash}p{.24\linewidth}p{.66\linewidth}}
\toprule Term & Meaning used in this manuscript\\\midrule\endhead
Local stability & Decay of small perturbations of a specified equilibrium, with no assertion of a global basin.\\
Tracking & The selected height error decays or satisfies an explicit disturbance envelope.\\
Internal dynamics & The remaining falling-speed equation after height and its derivative are assigned.\\
Invariant region & A set the constructed trajectory cannot leave under the specified law and hypotheses.\\
Force feasibility & The computed input and resulting pickup tension satisfy the given inequalities.\\
Global deterministic solution & A unique model trajectory on every forward time, within the invariant region and ideal sustained supply.\\
Maximal local SDE solution & A process defined until a positive lifetime, with an exhaustion-based explosion alternative and no claim of infinite lifetime.\\
Almost everywhere & Exceptional time sets may have zero Lebesgue measure.\\
Almost surely & Exceptional sample paths may have probability zero.\\
Stopped estimate & A bound evaluated before a random exit and a finite deterministic horizon.\\
Kernel-checked statement & A formal proposition proved from its displayed hypotheses in the pinned Lean environment.\\
Model validity & Whether the equations represent the apparatus, as determined by physical evidence outside the formal proof.\\\bottomrule
\end{longtable}

\FloatBarrier
\section{Counterexamples delimiting the guarantees}\label{app:counter}
The following examples identify concrete errors that arise if the scope conditions are omitted. The code also reproduces the first two examples, while the third follows directly from an analytical boundary check.

\subsection{A crown margin does not preserve forward transport}
With $H=2$, $g=9.81$, $\beta=0.02$, $h=0.1$, $v_1=10$, $w=0$, and $v_{\min}=3.2$, one has $B_C=1.8>0$, while
\begin{equation}
 f_2=9.81-\frac{25}{2.1}=-2.0947619\ldots<0.
\end{equation}
The input does not enter this boundary derivative. The crown inequality alone therefore gives no inward condition for the falling-speed boundary. The script starts instead at the interior state $(0.1,10,0.001)$ and uses the earlier unweighted nominal controller with constant reference $0.1$, gains $3.5$ and $8$, and the nominal crown filter with $\gamma=10$. The trajectory reaches $w=0$ at approximately $0.00048317\,\mathrm s$ while $B_C$ remains approximately $1.79182$. Integration stops at that event because continuation with negative $w$ would leave the assumed transport regime.

\subsection{The robust interval can be empty}
Take $H=2$, $h=3$, $v_1=6.3$, $w=0.1$, $\beta=0.02$, $v_{\min}=3.2$, and $\alpha\in[0.06,0.18]$. Both velocities and the height are positive and $B_C=0$, but
\begin{equation}
 L_B=-32.5458,\qquad U_B=-33.21248.
\end{equation}
The lower bound exceeds the upper bound. Although this example lies outside the conservative operating region used in Section~\ref{sec:numerics}, it disproves universal feasibility over the whole positive domain. The correct general statement is the conditional result in Proposition~\ref{prop:barrier}.

\subsection{Arbitrary bounded reference rates need not be feasible}
On the exact tracking manifold, $q=\dot r$. A smooth positive reference can have bounded but large derivative. For example $r(t)=0.2+0.1\sin(100t)$ has $\dot r(0)=10$, so at $w=0$ and $h=r(0)$,
\begin{equation}
 \dot w=9.81-\frac{100}{2.2}<0.
\end{equation}
Thus positivity and boundedness of the reference and its derivatives are insufficient on their own. Their numerical bounds must enter a feasibility condition such as \eqref{eq:lowerface}. The constant-reference Riccati result remains valid.

\FloatBarrier\section{Theorem-to-Lean correspondence and verification scope}\label{app:formal}
The current complete entry point is \lean{ChainFountainCompletion}. The v0.19.0 release includes the original mathematical core, the zero-gravity extension and the completion results listed below. Its combined compilation and transitive axiom audit use Lean~4.33.0-rc1 with the same pinned mathlib, MathFin and Brownian-motion dependencies as the baseline.

The source inventory contains 722 named authored theorems. Every one of the 1672 project-owned kernel declarations, including generated helpers, was included in the transitive axiom audit. The dependencies use only \texttt{propext}, \texttt{Classical.choice}, and \texttt{Quot.sound}, with no \texttt{sorry} or custom mathematical axiom in the project. These counts describe formal declarations, not independent contributions. The exact types, source hashes, dependency edges and audit result are supplied in \path{verification/completion/}.

The v0.16.0 source and its 557-theorem audit remain historical evidence. The later v0.18.0 audit checked 78 additional zero-gravity theorems. The present combined audit supersedes their separate-build limitation without altering their mathematical proof bodies.

\subsection{What has and has not changed in this edition}
The expanded preprint reorganizes and explains the existing mathematical results; the present completion also supplies formal proofs for the remaining explicit deductions. The ten numbered deterministic result statements use the original text and labels from the v0.16 manuscript. Three previously unnumbered stochastic results are now stated separately in Section~\ref{sec:stochastic}. The baseline mathematical statements and analytical figures are retained. The finite-source and geometric-obstruction deductions now have separately stated results, and the additional correspondence appears below. Appendix~\ref{app:zg-math} describes the zero-gravity development and its coverage within the combined audit. Assessing whether the prose matches the formal hypotheses still requires mathematical interpretation.

The v0.16 register has 66 claim groups: 60 checked groups, four numerical groups, one model group and one scope group. All mathematical obligations in that register are closed. The register covers the specified claims, leaving further extensions, explanatory prose and experimental adequacy outside that assertion of completeness. In particular, deterministic global continuation is proved under the operating-region assumptions, whereas the stochastic mechanical process has a maximal local lifetime which is not proved infinite.

The full research package includes an explicit crosswalk from this edition's labels to the v0.16 claim IDs, the unmodified formal project archive, its verification summary, and the original interactive dependency map. The map contains the shorter manuscript's source snapshot and line numbers, so its source coordinates refer to that manuscript. The crosswalk in this preprint supplies the section and page numbers for the expanded text.

The public v0.18 \href{https://github.com/luchesart/Chain-fountain/blob/main/Chain_Fountain_Lean_GitHub.zip}{Lean source archive} provides the preceding pinned project and zero-gravity extension. The v0.19 completion accompanies this edition as a separate updated source package. Its separate Mathlib-only browser companion contains 350 theorems, including the 78 zero-gravity results, and excludes the mechanical SDE, actual stochastic integrals, LIL and the new completion modules. It is retained as a historical browser companion. The \href{https://luchesart.github.io/Chain-fountain/Chain_Fountain_Zero_Gravity_Proof_Map.html}{interactive proof map} and \href{https://luchesart.github.io/Chain-fountain/Chain_Fountain_Dependency_Overview.pdf}{printable dependency overview} provide the declaration-level correspondence and dependencies.

\subsection{Entry points for the numbered results}
In the table, names without an explicit prefix belong to \texttt{ChainFountain}, while names prefixed by \texttt{S.} belong to \texttt{ChainFountainStochastic}. Long identifiers are broken only for typesetting. The supplied project and crosswalk give the full source, exact hypotheses and supporting declarations.
\begin{longtable}{>{\raggedright\arraybackslash}p{.27\linewidth}>{\raggedright\arraybackslash}p{.60\linewidth}}
\toprule Result and register group & Selected checked entry points\\\midrule\endhead
Energy balance (\ref{prop:energy}), P31 & \lean{energy_balance_ae}\\[4pt]
Steady-state input map (\ref{prop:setpoint}), P32 & \lean{steady_state_iff}\newline \lean{steady_actuator_interval}\\[4pt]
Natural local stability (\ref{thm:spectrum}), T41 & \lean{open_loop_local_exponential_unique}\newline \lean{transfer_physical_zero}\\[4pt]
Nominal error decay (\ref{prop:tracking}), P51 & \lean{closed_loop_error_equations}\newline \lean{nominal_error_estimate}\\[4pt]
Physical tracking region (\ref{thm:tube}), T61 & \lean{physical_region_invariant}\newline \lean{GlobalTrackingData.exists_unique_global_tracking}\\[4pt]
Full-state convergence (\ref{cor:fullstate}), C62 & \lean{GlobalTrackingData.exists_global_convergent}\\[4pt]
Force feasibility (\ref{prop:bounds}), P63 & \lean{robust_force_feasibility}\\[4pt]
Bounded uncertainty (\ref{thm:robust}), T71 & \lean{robust_error_estimate_ae}\newline \lean{GlobalTrackingData.exists_global_aemeasurable_admissible}\\[4pt]
Global construction and uniqueness (\ref{prop:global}), P72 & \lean{GlobalTrackingData.exists_unique_global_caratheodory}\newline \lean{GlobalTrackingData.exists_unique_classical_caratheodory}\\[4pt]
Barrier compatibility (\ref{prop:barrier}), P81 & \lean{barrier_feasible_iff}\newline \lean{filtered_mechanical_crown_ac}\\[4pt]
Constructed stopped energy (\ref{prop:stochastic-energy}), S116 & \lean{S.CutoffExhaustion.Family.maximal_stopped_energy_expectation}\\[4pt]
Actual integral freezing (\ref{thm:freezing}), S121 & \lean{S.LocalFreezingCompletion.local_freezing}\newline \lean{S.LocalFreezingCompletion.local_freezing_stoppedProcess}\\[4pt]
Mechanical crown LIL (\ref{thm:mechanical-lil}), S123 & \lean{S.MechanicalSmallTime.maximal_crown_normalized_lil}\newline \lean{S.MechanicalSmallTime.boundary_negative}\newline \lean{S.MechanicalSmallTime.exists_mechanical_crown_lil}\\[4pt]
\bottomrule\end{longtable}

The global wrappers take the reference, initial conditions and inequalities as their data and construct the solution from them. The AC/integral bridge and uniqueness theorems also remove the need to assume a state box for competing deterministic solutions. On the stochastic side, the construction establishes integral versions, cutoff compatibility, genuine exits and the maximal process before applying the energy and LIL conclusions.

For Theorem~\ref{thm:freezing}, \lean{S.LocalFreezingCompletion.local_freezing} proves the actual localized-integral conclusion without a global coefficient bound. The predictable and continuous-adapted integrals are identified by a proved simultaneous-in-time equality, and the stopping wrapper proves invariance of the small-time limit before any almost surely positive random stop. The filtration and continuous-version hypotheses are stated explicitly in Section~\ref{sec:stochastic}. The correspondence for Theorem~\ref{thm:mechanical-lil} still concerns the constructed maximal mechanical solution and its formal interface; unrestricted weak-solution uniqueness is not asserted.

\subsection{Completion of the remaining article deductions}
Identifiers below are prefixed by \lean{ChainFountain.ZeroGravity}, except for \lean{S.}, which denotes \lean{ChainFountainStochastic}.
\begin{longtable}{>{\raggedright\arraybackslash}p{.28\linewidth}>{\raggedright\arraybackslash}p{.61\linewidth}}
\toprule Article deduction & Checked entry points and precise scope\\\midrule\endhead
Finite source (\ref{cor:zg-depletion}) & \lean{constant_height_finite_supply}; \lean{remainingSupply_eq_zero_iff}. Actual integrated speed, unique depletion time and original mechanical equations while supply is positive.\\[5pt]
Closed-loop energy (\ref{eq:zg-loop-energy}) & \lean{loop_kinetic_integral}; \lean{loop_spatial_power_balance}; \lean{AECirculatingLoop.energy_balance}. Actual kinetic and force-power integrals, zero internal work, locally AC energy and a genuine integral power balance.\\[5pt]
Angular momentum and reaction torque & \lean{planar_material_torque_balance_C1}; \lean{planar_tangential_force_resultant_zero}. Actual material integrals with proved integrability, periodic phase removal and the time derivative of angular momentum.\\[5pt]
Shape obstruction (\ref{prop:zg-shape-obstruction}) & \lean{circulation_shape_attracted_iff}; \lean{circulation_speed_shape_obstruction}; \lean{arbitrarily_close_nonattracting_circulation}. Hausdorff distance of compact spatial images in a fixed frame, equal speed outputs and arbitrarily close translated examples.\\[5pt]
Arclength coordinates (Appendix~\ref{app:zg-math}) & \lean{closed_regular_curve_arclength}. A global $C^2$ unit-speed parameterization constructed from the integral of speed, with the same range and period equal to integrated length.\\[5pt]
Localized freezing (\ref{thm:freezing}) & \lean{S.LocalFreezingCompletion.local_freezing}; \lean{S.LocalFreezingCompletion.local_integral_lil}; \lean{S.LocalFreezingCompletion.local_freezing_stoppedProcess}. Actual stochastic integrals, no global coefficient bound, natural-filtration continuous versions and a positive initial stopping interval.\\
\bottomrule\end{longtable}

The accompanying \path{Chain_Fountain_Completion_Map.html} shows these claim groups, exact theorem types, source proofs and dependencies extracted from the elaborated terms. Its verification status requires both a successful combined audit and matching source hashes. The physical model, numerical experiments and future spatial-stability questions are outside this formal mathematical coverage.

\subsection{How to read a dependency graph}
The explanatory graphs in Figures~\ref{fig:proofroute} and~\ref{fig:stochasticroute} show the main mathematical route. The separate interactive map has theorem edges extracted from proof dependencies and module summaries. Importing a module does not mean that every theorem in it is used. A proof may depend on a completed part of a library that also contains unrelated unfinished work, so the transitive axiom audit follows the dependencies of the theorem itself.

Matching a paper claim to a formal declaration requires agreement in four respects: the quantified objects, all positivity, boundedness and regularity hypotheses, the time interval and solution class, and the conclusion with any stopping or almost-sure qualification. These comparisons establish the mathematical correspondence that a shared name alone cannot provide. The crosswalk records the stated limits alongside each claim.

\subsection{A reading map for the mathematical tools}
\begingroup\small
\begin{longtable}{>{\raggedright\arraybackslash}p{.24\linewidth}p{.37\linewidth}p{.29\linewidth}}
\toprule Paper task & Essential mathematical step & Background in this preprint\\\midrule\endhead
Energy balance & Product rule for moving masses and AC chain rule & Sections~\ref{sec:derive}, \ref{app:ac}\\
Natural stability & Cubic root exclusion and quadratic nonlinear remainder bound & Appendix~\ref{app:hurwitz}\\
Error cancellation & Invertible coordinates and mixed-term cancellation & Section~\ref{sec:tracking}\\
Error envelopes & Integrating factors and weighted norms & Appendix~\ref{app:comparison}\\
Internal-state bounds & Face factorization and weak comparison & Appendix~\ref{app:faces}\\
Global deterministic result & Local construction, compact continuation, unrestricted uniqueness & Appendices~\ref{app:ac}, \ref{app:existence}\\
Force/reference feasibility & Interval bounds and derivative maxima & Appendix~\ref{app:constraints}\\
Robust estimate & Young inequality, comparison, finite-energy integration & Section~\ref{sec:robust}, Appendix~\ref{app:comparison}\\
Crown filter & Intersection and projection of scalar input constraints & Appendix~\ref{app:constraints}\\
Stopped stochastic energy & It\^o formula, integrability, stopping & Appendices~\ref{app:ito}, \ref{app:generator}\\
Actual-integral freezing & Exponential maximum, Borel--Cantelli, prefix locality & Appendix~\ref{app:freezing}\\
Mechanical boundary exit & Sharp Brownian LIL, coefficient freezing, drift negligibility & Appendix~\ref{app:lil}\\\bottomrule
\end{longtable}
\endgroup

\FloatBarrier\section{Reproduction and claim map}
The public \href{https://github.com/luchesart/Chain-fountain}{project repository} hosts the Lean source archive and the browser companion. The minimal arXiv source package contains the complete manuscript, its TeX inputs and the figures needed for compilation. Compile \texttt{main.tex} with pdfLaTeX until the references settle. The separate full research package also contains Python source, dependencies, CSV trajectories, per-realization metrics and environment records. In that package, \texttt{python code/reproduce.py} regenerates the numerical plots and tables with NumPy, SciPy and Matplotlib. Integration settings are recorded in \texttt{results/metrics.json}.

\paragraph{Interactive companion.}
\href{https://luchesart.github.io/Chain-fountain/}{\textit{Chain Fountain Lab}} is publicly available at \url{https://luchesart.github.io/Chain-fountain/}, with English and Bulgarian controls. It integrates the three-state equations, animates material transport, and compares the selected feedback law with an unactuated trajectory from the same initial state. Readers can vary the reference, gains, bounded pickup fluctuations, force limits, actuator lag and finite chain supply. The displayed sufficient inequalities are evaluated independently of the numerical trajectory.

The links and crown shape are schematic transport markers, not a discrete-link mechanical simulation. The random-input option supplies a bounded ODE input and does not simulate the Brownian model of Section~\ref{sec:stochastic}. Adaptive RK4 trajectories were checked against an independently implemented DOP853 reference, the exact nominal error system and step refinement. These numerical checks do not amount to kernel verification of the browser solver. The hosted page runs entirely in the browser. Source files and test records are included in the separate companion package.

The zero-gravity workspace adds an open-transport mode, an exact continuum circulation family, and a separate elastic bead-ring experiment. Section~\ref{sec:zg-simulation} explains the assumptions behind each construction and its role in examining the dynamics.

\begin{table}[htbp]
\centering\small
\begin{tabular}{p{.25\linewidth}p{.64\linewidth}}
\toprule
Result & Exact scope \\
\midrule
Energy and setpoint map & Identities within the vertical reduced model; no new microscopic mechanism asserted.\\
Local stability & Uncontrolled positive three-state equilibrium, stated coefficient range.\\
Nominal tracking & Exponential error decay; full physical operating conditions supplied by the invariant region and force bounds.\\
Bounded uncertainty & Weighted disturbance estimate and robust invariant region, for unfiltered force feedback under the stated inequalities.\\
Global existence & Unique global solution from initial data in the invariant region; classical and measurable existence/physical bounds Lean-checked; unrestricted AC uniqueness Lean-checked in both nominal and robust branches.\\
Crown filter & Exact pointwise feasibility and conditional crown/pickup guarantees; no automatic full-domain invariance.\\
White noise & Constructed maximal local solution, stopped estimate and actual-integral crown LIL; almost-sure boundary invariance is ruled out when the normal diffusion is nonzero.\\
Numerics & Reproducible reduced-model illustrations; no experimental validation or universal performance ranking.\\
\bottomrule
\end{tabular}
\caption{Map of established statements and their limits.}
\end{table}

\FloatBarrier
\FloatBarrier\section{Mathematics and formal scope of the zero-gravity extension}\label{app:zg-math}
\subsection{Material coordinates and the chain rule}
The label $a$ follows a material point of the string as it moves through space. Periodicity identifies $a$ and $a+L$, while inextensibility requires the material arclength metric to satisfy $|X_a|=1$. For $X=R(a+\theta)$, $R'$ is the unit tangent and $R''$ is the curvature vector. Differentiation at fixed material label produces both the normal acceleration $c^2R''$ and the tangential acceleration $\dot cR'$. The normal term accounts for the change of direction along the curve. Omitting it would incorrectly predict zero required tension at constant speed.

An arbitrary smooth regular closed curve $Q(b)$ can be parameterized by arclength: define $\ell(b)=\int_0^b|Q'(\xi)|\,\dd\xi$ and $R(a)=Q(\ell^{-1}(a))$. In a renderer this inverse may be approximated by a sufficiently fine cumulative-length table. That numerical approximation is distinct from the exact unit-speed hypothesis in Proposition~\ref{prop:zg-loop}.

\subsection{The integral solution and its estimate}
For a locally integrable disturbance, the integrating-factor identity is
\[
 e^{\mu t}e(t)=e(0)+\int_0^t e^{\mu\tau}d(\tau)\,\dd\tau.
\]
Absolute continuity allows the derivative to exist almost everywhere while ensuring that the fundamental theorem of calculus recovers increments by integration. With $|d|\le D$ almost everywhere, the triangle inequality for the integral gives the exponential envelope in \eqref{eq:zg-bound}. This formula explicitly constructs a solution of the scalar equation for all forward time, so global existence follows from the construction itself.

\paragraph{From measurable data to a solution.}
On every $[0,T]$, almost-everywhere measurability and the essential bound imply integrability of $d$. The continuous integrating factor is bounded on that compact interval, so $s\mapsto e^{\mu s}d(s)$ is integrable there as well. Its indefinite integral is absolutely continuous, with derivative equal to the integrand almost everywhere. Multiplication by $e^{-\mu t}$ then constructs a locally absolutely continuous error $e$. Consequently, $c=r+e$ has initial value $c_0$ and satisfies the controlled equation almost everywhere. Removing the union of the exceptional null sets obtained on the integer intervals $[0,n]$ gives one full-measure set on the entire forward half-line. The assumptions concern only $t\ge0$, and the construction requires no extension of the disturbance before zero.

\paragraph{Why uniqueness includes every absolutely continuous solution.}
For any competing locally absolutely continuous solution, the product $e^{\mu t}e(t)$ is absolutely continuous on each compact interval and has derivative $e^{\mu t}d(t)$ almost everywhere. The reverse direction of the fundamental theorem of calculus therefore gives the displayed integrating-factor identity for every such solution. This establishes the integral representation without imposing it as an additional condition on the solution class. Equal initial values then give equal forward trajectories. Applying the integral triangle inequality to the representation yields
\[
 |e(t)|\le e^{-\mu t}\!\left(|e(0)|+D\int_0^t e^{\mu s}\,\dd s\right)
 =|e(0)|e^{-\mu t}+\frac D\mu(1-e^{-\mu t}).
\]
This estimate and all resulting speed, tension and command bounds hold at every $t\ge0$, even though the differential equation and the forcing bound are almost-everywhere statements.

\paragraph{From the scalar solution to continuum momentum.}
The continuous extension $\widetilde c(t)=c(\max\{0,t\})$ defines
$\theta(t)=\theta_0+\int_0^t\widetilde c(s)\,\dd s$ with a classical derivative at every time. For a $C^2$ curve, the material velocity
$V(a,t)=c(t)R'(a+\theta(t))$ is locally absolutely continuous in forward time. On one full-measure set, independent of $a$, its derivative equals
\[
 V_t(a,t)=(u(t)-\gamma c(t)+d(t))R'(a+\theta(t))
             +c(t)^2R''(a+\theta(t)).
\]
The spatial derivative of $T X_a$ is the second term multiplied by $\lambda$. The temporal balance therefore holds for all material labels simultaneously outside that single exceptional set. The coordinatewise fundamental theorem of calculus also gives the genuine vector integral momentum identity
\[
 \lambda\bigl(V(a,t)-V(a,0)\bigr)=
 \int_0^t\left[\partial_a(TX_a)(a,s)
     +\lambda\bigl(u(s)-\gamma c(s)+d(s)\bigr)X_a(a,s)\right]\,\dd s.
\]
The formal three-coordinate integral identity uses the usual coordinate representation of vectors. Physical unit speed and inextensibility are stated separately in a real inner-product space, in particular Euclidean two- or three-space.

Strict positivity and exact reconstruction place the constructed motion in the physical region. This argument makes no assumption that other spatial solutions approach the invariant family. In the general PDE, tension is a Lagrange multiplier enforcing inextensibility, so it cannot be chosen independently of $X$. The expression $\lambda c^2$ is justified only on the displayed family.

\subsection{What is checked, and what is not}
The Lean module register below distinguishes the earlier classical-solution entry points from the extension to measurable inputs and absolutely continuous solutions. It identifies the checked declarations for the scalar construction, uniqueness, integral representation, physical bounds and continuum reconstruction. This formal development starts from the classical circulation family described by Healey and Papadopoulos~\cite{healey1990} and uses the standard integrating-factor solution of a scalar linear equation. Its purpose is to verify how those ingredients yield the stated controlled motion and its physical bounds, with every regularity and solution-class assumption explicit. Each formal proof establishes the corresponding deduction from its hypotheses. The physical adequacy of the hypotheses and the attribution of the underlying mechanics require the separate comparisons made in Section~\ref{sec:zerog}. The completion modules now prove finite-supply depletion, the planar angular-momentum and torque integrals with their integrability, kinetic-energy and power balances, geometric non-attraction in Hausdorff distance, and global arclength reparameterization. Their entry points are listed in Appendix~\ref{app:formal} and included in the common audit.
The v18 verification compiled all 5 zero-gravity modules and their local import closure, checking \textbf{78 authored theorems}: 39 results from v17 and 39 theorems added in v18. The imported absolute-continuity bridge was compiled in the same verification. Every target-owned kernel declaration passed a transitive axiom audit allowing only \texttt{propext}, \texttt{Classical.choice} and \texttt{Quot.sound}. The independent target is \texttt{lake build ChainFountainZeroGravity}. The archived 557-theorem baseline audit provides the earlier verification evidence. The full stochastic baseline was not rebuilt for v18; the present v0.19 combined audit now includes it together with the completion modules. The table below retains the original ZG-O, ZG-C and ZG-A correspondence; the new deductions are listed in Appendix~\ref{app:formal}.

\medskip

\textbf{Measurable-input coverage.} The ZG-A entries formalize the bounded-measurable/locally-absolutely-continuous statements of Proposition~\ref{prop:zg-loop} and Theorem~\ref{thm:zg-tracking}. The construction gives a scalar trajectory and its phase. An integrating-factor formula derived for every AC competitor establishes forward uniqueness and pointwise physical bounds. The reconstructed continuum balance holds on one full-measure time set shared by all material labels. A separate three-coordinate result proves the vector integral momentum identity. The earlier ZG-C entries have their original classical hypotheses. The ZG-A entries state the measurable-input extensions separately.

\medskip

All identifiers below have namespace \lean{ChainFountain.ZeroGravity}. The interactive map provides exact hypotheses, source text, authored declarations and generated helpers, and dependencies extracted from the elaborated Lean terms.

\begin{longtable}{>{\raggedright\arraybackslash}p{.30\linewidth}>{\raggedright\arraybackslash}p{.60\linewidth}}

\toprule Paper claim & Selected Lean entry points and scope\\\midrule\endhead

ZG-O1 (\ref{thm:zg-budget}): Positive-speed lower face is impossible at zero gravity & \lean{lower_face_impossible}\\[5pt]

ZG-O2 (\ref{thm:zg-budget}): Exact falling momentum derivative and crown dissipation & \lean{falling_momentum_hasDerivAt}\newline \lean{falling_momentum_rate_le}\\[5pt]

ZG-O3 (\ref{thm:zg-budget}): Finite positive-margin horizon and global impossibility & \lean{positive_margin_time_bound}\newline \lean{global_positive_margin_impossible}\\[5pt]

ZG-O4 (\ref{prop:zg-decay}): No positive-flow equilibrium & \lean{no_positive_flow_equilibrium}\\[5pt]

ZG-O5 (\ref{prop:zg-decay}): Explicit constant-height rational decay solves original plant & \lean{decaySpeed_pos}\newline \lean{constant_height_solution}\\[5pt]

ZG-O6 (\ref{prop:zg-decay}): Zero-gravity energy identity & \lean{energy_balance_zeroGravity}\\[5pt]

ZG-C1 (\ref{thm:zg-tracking}): Exact scalar closed-loop error equation & \lean{circulation_error_hasDerivAt}\\[5pt]

ZG-C2 (\ref{thm:zg-tracking}): Exact nominal scalar evolution & \lean{circulation_nominal_global_exists}\\[5pt]

ZG-C3 (\ref{thm:zg-tracking}): Constructive global continuous-forcing solution and uniqueness & \lean{circulation_forced_global_exists}\newline \lean{circulation_scalar_unique}\newline Continuous forcing. Global scalar construction and finite-interval uniqueness.\\[5pt]

ZG-C4 (\ref{thm:zg-tracking}): Robust exponential error envelope and positive circulation & \lean{circulation_error_bound}\newline \lean{circulation_positive}\newline Classical or one-sided derivatives. See ZG-A for the a.e. result.\\[5pt]

ZG-C5 (\ref{thm:zg-tracking}): Tension and actuator bounds & \lean{circulation_command_bound}\newline \lean{circulation_tension_upper}\newline \lean{loop_tension_positive}\\[5pt]

ZG-C6 (\ref{prop:zg-loop}): Actual continuum derivative reconstruction & \lean{loop_continuum_balance}\newline Classical derivatives. See ZG-A for the constructed AC motion.\\[5pt]

ZG-C7 (\ref{prop:zg-loop}): Periodicity and inextensibility & \lean{loop_inextensible}\newline \lean{loop_periodic}\newline Periodicity and unit speed are explicit curve hypotheses.\\[5pt]

ZG-A1 (\ref{thm:zg-tracking}): Constructive global-forward measurable/AC circulation & \lean{circulation_forward_measurable_global_exists}\newline \lean{circulationACSolution_integral}\\[5pt]

ZG-A2 (\ref{thm:zg-tracking}): AC representation and unrestricted forward uniqueness & \lean{circulation_ac_closedLoop_unique}\newline \lean{circulation_ac_variation_of_constants}\\[5pt]

ZG-A3 (\ref{thm:zg-tracking}): Sharp a.e. disturbance estimate and pointwise physical bounds & \lean{circulation_measurable_global_wellposed}\newline \lean{circulation_ac_error_bound}\\[5pt]

ZG-A4 (\ref{prop:zg-loop}): Constructed C1 phase and absolutely continuous material velocity & \lean{circulationPhase_contDiff}\newline \lean{loop_material_velocity_ac}\\[5pt]

ZG-A5 (\ref{prop:zg-loop}): Actual continuum PDE on one common full-measure time set & \lean{loop_continuum_balance_ae}\\[5pt]

ZG-A6 (\ref{prop:zg-loop}): Genuine coordinate integral momentum balance & \lean{loop_ac_integral_balance}\\[5pt]

ZG-A7 (\ref{prop:zg-loop}): Unconditional construction of the measurable-driven closed loop & \lean{measurable_closed_loop_exists}\\[5pt]

\bottomrule\end{longtable}

The graphical companion shows the dependencies between the new declarations alongside their source text and the corresponding paper labels. A separate historical v0.16 map preserves the theorem identifiers and manuscript coordinates from that version. Neither map presents the elastic bead simulator as a formalization of the inextensible PDE. Further work is needed on general spatial well-posedness and transverse stability, on localized physical actuator and contact dynamics, and on a closed-loop spatial SDE describing stochastic shape evolution.

\section*{AI assistance and verification status}
An initial draft was developed with Google Gemini. ChatGPT assisted with auditing derivations, the invariant-region conditions, numerical implementation, formal proofs, and preparation of this expanded exposition and its figures with multiple iterations under the guidance of the author. The author has checked the proofs and prepared an \href{https://luchesart.github.io/Chain-fountain/Chain_Fountain_Proof_Map.html}{interactive map of the formalization}. The map connects statements and proofs in the paper with the corresponding Lean declarations, their explicit hypotheses, and the dependencies between theorems and lemmas.

The accompanying \href{https://luchesart.github.io/Chain-fountain/}{interactive simulator} allows readers to explore the model and the effect of the controller. The \href{https://github.com/luchesart/Chain-fountain/blob/main/Chain_Fountain_Lean_GitHub.zip}{Lean source archive} provides the formal project and its zero-gravity extension. A \href{https://luchesart.github.io/Chain-fountain/Chain_Fountain_Dependency_Overview.pdf}{printable dependency overview} complements the interactive map.

This edition retains the existing numerical experiment outputs without rerunning the experiments. The formal completion and combined Lean audit are recorded in Appendix~\ref{app:formal}. The new analytical figures are generated from the displayed formulas by our scripts. The updated map organizes the recorded verification results and does not constitute a new numerical validation.

The v0.16.0 source and its 557-theorem audit remain available as historical evidence. Appendix~\ref{app:zg-math} presents the zero-gravity correspondence and its exact scope within the combined audit. Baseline correspondence remains in Appendix~\ref{app:formal}. Lean checks deductions from explicit hypotheses. Model adequacy, numerical accuracy and empirical validity each require separate evidence. The previously checked variable-diffusion transfer is preserved. Final author review of the text, correspondence, attribution and implementation remains necessary before submission. Experimental validation is incomplete.

\end{document}